\documentclass[12pt]{article}
\usepackage{latexsym,amssymb,amsmath}
\begin{document}

\baselineskip=20pt

\newcommand{\la}{\langle}
\newcommand{\ra}{\rangle}
\newcommand{\psp}{\vspace{0.4cm}}
\newcommand{\pse}{\vspace{0.2cm}}
\newcommand{\ptl}{\partial}
\newcommand{\dlt}{\delta}
\newcommand{\sgm}{\sigma}
\newcommand{\al}{\alpha}
\newcommand{\be}{\beta}
\newcommand{\G}{\Gamma}
\newcommand{\gm}{\gamma}
\newcommand{\vs}{\varsigma}
\newcommand{\lmd}{\lambda}
\newcommand{\td}{\tilde}
\newcommand{\vf}{\varphi}
\newcommand{\rd}{\mbox{Rad}}
\newcommand{\ad}{\mbox{ad}}
\newcommand{\stl}{\stackrel}
\newcommand{\ol}{\overline}
\newcommand{\ul}{\underline}
\newcommand{\es}{\epsilon}
\newcommand{\dmd}{\diamond}
\newcommand{\clt}{\clubsuit}
\newcommand{\vt}{\vartheta}
\newcommand{\ves}{\varepsilon}
\newcommand{\dg}{\dagger}
\newcommand{\kn}{\mbox{ker}}
\newcommand{\for}{\mbox{for}}
\newcommand{\rad}{\mbox{Rad}}

\begin{center}{\Large \bf New Generalized Simple Lie Algebras of Cartan}\end{center}
\begin{center}{\Large \bf  Type over a Field with Characteristic 0}
\end{center}
\vspace{0.2cm}

\begin{center}{\large Xiaoping Xu}\end{center}
\begin{center}{Department of Mathematics, The Hong Kong University of Science \& Technology}\end{center}
\begin{center}{Clear Water Bay, Kowloon, Hong Kong}\footnote{Research supported
 by Hong Kong RGC Competitive Earmarked Research Grant HKUST709/96P.}\end{center}

\vspace{0.3cm}

\begin{center}{\Large \bf Abstract}\end{center}
\vspace{0.2cm}

{\small We construct four new series of generalized simple Lie algebras of Cartan type, using the mixtures of grading operators and down-grading operators. Our results in this paper are further generalizations of those in Osborn's work ``New simple infinite-dimensional Lie algebras of characteristic 0.''}

\section{Introduction}

The four well-known series of simple Lie algebras of Cartan type were constructed from the derivation algebra of the polynomial algebra in several variables and its subalgebras preserving certain differential forms. The abstract definitions by derivations of generalized Lie algebras of  Cartan type appeared in Kac's work [Ka1]. However, it is still a question of how to construct new explicit generalized simple Lie algebras of Cartan type.
Kawamoto [K] introduced new simple Lie algebras of generalized Witt type by changing the polynomial algebra in several variables to the group algebra of an additive subgroup of $\Bbb{F}^n$  and considering the Lie subalgebra of its derivation algebras generated by the grading operators, where $n$ is a positive integer and $\Bbb{F}$ is a field with characteristic 0. One can view the operators of taking partial derivatives of  the polynomial algebra in several variables as {\it down-grading operators}. Starting from the derivation subalgebra  generated by the grading operators and down-grading operators of the tensor algebra of the group algebra of the direct sum of finite number of additive subgroups of $\Bbb{F}$ with the polynomial algebra in several variables, Osborn [O2] constructed new four series of  generalized simple Lie algebras of Cartan type. In [DZ], the authors generalized Kawamoto's work by picking out certain subalgebras, whose typical examples are the derivation subalgebra generated by the grading operators and down-grading operators of the tensor algebra of the group algebra of an additive subgroup of $\Bbb{F}^n$ with the polynomial algebra in several variables. The work in [DZ] was also a generalization of Osborn's generalized Witt algebra in [O2]. Passman [P] gave a certain necessary and sufficient condition on derivations for a Lie algebra of generalized Witt type to be simple.

In [O1], Osborn gave a classification of infinite-dimensional simple Novikov algebras with an idempotent element over a field with characteristic 0. We observed in [X3] that the three classes of Osborn's classified simple Novikov algebras can be rewritten in one form by using the sum of a grading operator and a down-grading operator. In other words, these three classes of Novikov algebras can be viewed as one class in terms of our notions. This observation leaded us to construct in [X3] a much larger class of simple Novikov algebras, including  Osborn's classified simple Novikov algebras as very special cases. We also observed that the commutator Lie algebras of Osborn's classified simple Novikov algebras are rank-one simple Lie algebras of generalized Witt type. In fact, by considering the Lie algebras induced by the Hamiltonian operators corresponding to the simple Novikov algebras with an idempotent element, we found in [X3] a new family of infinite-dimensional Lie algebras. Moreover, a new family of infinite-dimensional Lie superalgebras were constructed based on the Hamiltonian superoperators (cf. [X2]) corresponding to the Novikov-Poisson algebras (cf. [X1]) whose Novikov algebraic structures are those classified in [O1]. Recently, we classified in [X4] quadratic conformal superalgebras by certain compatible pairs of a Lie superalgebra and a Novikov superalgebra. Six general constructions of such pairs  were given. Moreover, we classified  such pairs related to simple Novikov algebras with an idempotent element. As the consequences of this classification, new families of Lie algebras were found. One of the motivations of this paper is to understand the simplicity of these Lie algebras. Another motivation is to understand the simplicity of the Lie algebras generated by our constructed quadratic conformal superalgebras in [X4] (e.g., cf. Section 4.1 in [X5] for the constructions of these Lie algebras). 

Our main results in this paper are the constructions and proofs of four new series of generalized simple Lie algebras of Cartan type based on the derivation subalgebra generated by the mixtures of the grading and down-grading operators of the tensor algebra of the group algebra of an additive subgroup of $\Bbb{F}^n$ with the polynomial algebra in several variables. Our Lie algebras of type H contain some non-derivation ingredients. 

For the convenience of the reader's understanding this work, below we shall present the constructions of the four series of simple Lie algebras of Cartan type. 

Throughout this paper, let $\Bbb{F}$ be a field  with characteristic 0. All the vector spaces are assumed over $\Bbb{F}$. Denote by $\Bbb{Z}$ the ring of integers and by $\Bbb{N}$ the set of natural numbers $\{0,1,2,3,...\}$. We shall always identify $\Bbb{Z}$ with $\Bbb{Z}1_{\Bbb{F}}$ when the context is clear.

Let ${\cal  A}=\Bbb{F}[t_1,t_2,...,t_n]$ be the algebra of polynomials in $n$ variables. A {\it derivation} $\ptl$  of ${\cal A}$ is linear transformation of ${\cal A}$ such that
$$\ptl(uv)=\ptl(u)v+u\ptl(v)\qquad\for\;\;u,v\in{\cal A}.\eqno(1.1)$$
The typical derivations are $\{\ptl_{t_1},\ptl_{t_2},...,\ptl_{t_n}\}$,  the operators of taking partial derivatives. The space $\mbox{Der}\:{\cal A}$ of all the derivations of ${\cal A}$ forms a Lie algebra. Identifying the elements of 
${\cal A}$ with their corresponding multiplication operators, we have
$$\mbox{Der}\:{\cal A}=\sum_{i=1}^n{\cal A}\ptl_{t_i},\eqno(1.2)$$
which forms a simple Lie algebra. The Lie algebra $\mbox{Der}\:{\cal A}$ is called a {\it Witt algebra of rank} $n$, usually denoted as ${\cal W}(n,\Bbb{F})$. The Lie algebra ${\cal W}(n,\Bbb{F})$ acts on the Grassmann algebra $\hat{\cal A}$ of differential forms on ${\cal A}$ as follows.
$$\ptl(df)=d(\ptl(f)),\;\;\ptl(\omega\wedge \nu)=\ptl(\omega)\wedge\nu +\omega\wedge\ptl(\nu)\eqno(1.3)$$
for $f\in{\cal A},\;\omega,\nu\in\hat{\cal A},\;\ptl\in {\cal W}(n,\Bbb{F})$. Set
$${\cal S}(n,\Bbb{F})=\{\ptl\in {\cal W}(n,\Bbb{F})\mid \ptl(dt_1\wedge dt_2\wedge\cdots\wedge dt_n)=0\}.\eqno(1.4)$$
Assume that $n=2k$ is an even integer. Define
$${\cal H}(2k,\Bbb{F})=\{\ptl\in {\cal W}(n,\Bbb{F})\mid \ptl(\sum_{i=1}^kdt_i\wedge dt_{k+i})=0\}.\eqno(1.5)$$
Suppose that $n=2k+1$ is an odd integer. We let
\begin{eqnarray*}{\cal K}(2k+1,\Bbb{F})&=&\{\ptl\in {\cal W}(n,\Bbb{F})\mid \ptl(dt_{2k+1}+\sum_{i=1}^k(t_idt_{k+i}-t_{k+i}dt_i)\\& &\in {\cal A}(dt_{2k+1}+\sum_{i=1}^k(t_idt_{k+i}-t_{k+i}dt_i)\}.\hspace{5.3cm}(1.6)\end{eqnarray*}
The subspaces ${\cal S}(n,\Bbb{F}),\;{\cal H}(2k,\Bbb{F})$ and ${\cal K}(2k+1,\Bbb{F})$ form simple Lie subalgebras, which are called the Lie algebras of {\it Special type, Hamiltonian type} and {\it Contact type}, respectively. For convenience, we simply call the Lie algebras  ${\cal W}(n,\Bbb{F}),\;{\cal S}(n,\Bbb{F}),\;{\cal H}(2k,\Bbb{F})$ and ${\cal K}(2k+1,\Bbb{F})$ the Lie algebras of {\it type} W, S, H and K, respectively. It can be observed that the simplicity of the Lie algebras of type S, H and K is determined only by their concrete presentations of elements (e.g., cf. [SF]) and does not have any direct relations with their defining differential forms. This essentially gives us rooms to generalize these algebras.

We shall process our constructions and proofs section by section according to the Cartan type W, S, H and K. 

\section{Algebras of Type W}

In this section, we shall construct and prove a  new class of generalized simple Lie algebras of Witt type.

Given $m,n\in\Bbb{Z}$ with $m<n$, we shall use the following notion 
$$\ol{m,n}=\{m,m+1,m+2,...,n\}\eqno(2.1)$$
throughout this paper. We also treat $\ol{m,n}=\emptyset$ when $m>n$.
\psp

{\bf Definition 2.1}.  Let ${\cal A}$ be a commutative  associative algebra with an identity $1_{\cal A}$ and  let $\{\ptl_1,...,\ptl_n\}$ be $n$ linearly independent and mutually commutative derivations of ${\cal A}$ such that
$$\Bbb{F}1_{\cal A}+\sum_{p=1}^n\ptl_p({\cal A})={\cal A},\qquad \{u\in{\cal A}\mid \ptl_1(u)=\cdots =\ptl_n(u)=0\}=\Bbb{F}1_{\cal A}.\eqno(2.2)$$
Identify the elements of ${\cal A}$ with their corresponding multiplication operators. 
The following subspace of derivations
$${\cal W}=\sum_{i=1}^n{\cal A}\ptl_i\eqno(2.3)$$
 forms a Lie subalgebra of the Lie algebras of all the derivations of ${\cal A}$. We call ${\cal W}$ a {\it generalized Lie algebra of  Witt type}. In fact, its Lie bracket is given by
$$[\sum_{p=1}^nu_p\ptl_p,\sum_{q=1}^nv_q\ptl_q]=\sum_{p,q=1}^n(u_p\ptl_p(v_q)-\ptl_p(u_q)v_p)\ptl_q\eqno(2.4)$$
for $u_p,v_q\in{\cal A}$.
\psp

Next we shall give the specific construction of our generalized simple Lie algebras of Witt type. Let $n$ be a positive integer.  Pick
$${\cal J}_p\in\{\{0\},\Bbb{N}\}\qquad\for\;\;p\in\ol{1,n}.\eqno(2.5)$$
Let $\G$ be a torsion-free abelian group and let $\{\vf_p\mid p\in\ol{1,n}\}$ be $n$ additive group homomorphisms from $\G$ to $\Bbb{F}$ such that
$$\bigcap_{p\neq q\in\ol{1,n}}\kn_{\vf_q}\setminus\kn_{\vf_p}\neq\emptyset\qquad\mbox{if}\;\;{\cal J}_p=\{0\}\;\;\for\;\;p\in\ol{1,n},\eqno(2.6)$$
$$\bigcap_{p=1}^n\kn_{\vf_p}=\{0\}.\eqno(2.7)$$
Condition (2.7) implies that $\G$ is isomorphic to an additive subgroup of $\Bbb{F}^n$.

Set 
$$\vec{\cal J}={\cal J}_1\times{\cal J}_2\times\cdots\times{\cal J}_n,\eqno(2.8)$$
where the addition of $\vec{\cal J}$ is defined componentwise. Moreover, we denote: 
$$i_{[p]}=(0,..,\stackrel{p}{i},0,...,0)\qquad\for\;\;i\in{\cal J}_p.\eqno(2.9)$$
Let ${\cal A}$ be a vector space with a basis 
$$\{x^{\al,\vec{i}}\mid (\al,\vec{i})\in\G\times\vec{\cal J}\}.\eqno(2.10)$$
We define an algebraic operation $\cdot$ on ${\cal A}$ by
$$x^{\al,\vec{i}}\cdot x^{\be,\vec{j}}=x^{\al+\be,\vec{i}+\vec{j}}\qquad\for\;\;(\al,\vec{i}),(\be,\vec{j})\in\G\times\vec{\cal J}.\eqno(2.11)$$
Then ${\cal A}$ forms an associative algebra with an identity element $x^{0,\vec{0}}$, which will be simply denoted by $1$ in the rest of this paper.
Moreover, we define $\ptl_{\vf_p},\hat{\ptl}_p\in \mbox{End}\;{\cal A}$ for $p\in\ol{1,n}$ by:
$$\ptl_{\vf_p}(x^{\al,\vec{i}})=\vf_p(\al)x^{\al,\vec{i}},\;\;\hat{\ptl}_p(x^{\al,\vec{i}})=i_px^{\al,\vec{i}-1_{[p]}}\eqno(2.12)$$
for $(\al,\vec{i})\in\G\times\vec{\cal J}$ (cf. (2.9)), where we adopt the convention that if a notion is not defined but technically appears in an expression, we always treat it as zero; for instance, $x^{\al,-1_{[1]}}=0$ for any $\al\in\G$.  Then $\{\ptl_{\vf_p},\hat{\ptl}_q\mid p,q\in\ol{1,n}\}$ are derivations of ${\cal A}$. Furthemore,
$$\hat{\ptl}_p^{i_p+1}(x^{\al,\vec{i}})=0\qquad\for\;\;p\in\ol{1,n},\;(\al,\vec{i})\in\G\times\vec{\cal J}.\eqno(2.13)$$
Thus $\ptl_{\vf_p}$ is a grading operator and $\hat{\ptl}_p$ is a down-grading operator.
We set
 $$\ptl_p=\ptl_{\vf_p}+\hat{\ptl}_p\qquad\for\;\;p\in\ol{1,n}\eqno(2.14)$$
and define $({\cal W},[\cdot,\cdot])$ by (2.3) and (2.4). It can be verified that (2.2) holds.

\psp

{\bf Theorem 2.2}. {\it The Lie algebra} $({\cal W},[\cdot,\cdot])$ {\it is simple.}
\psp

{\it Proof}. It can be verified that our construction ingredients satisfy Passman's simplicity conditions of the Lie algebras of Witt type (cf. [P]).$\qquad\Box$
\psp

{\bf Example}. Let $n=n_1+n_2$ with $n_1,n_2\in\Bbb{N}$. We pick
$${\cal J}_p=\{0\},\;\;{\cal J}_{n_1+q}=\Bbb{N}\qquad\for\;\;p\in\ol{1,n_1},\;q\in\ol{1,n_2}\eqno(2.15)$$
(cf. (2.5)). Let $m\in\Bbb{N}$ be such that $n_1\leq m\leq n$ and define
$$\zeta_p(\al_1,\al_2,...,\al_m)=\al_p\qquad\for\;\;p\in\ol{1,m},\;(\al_1,\al_2,...,\al_m)\in\Bbb{F}^m.\eqno(2.16)$$
Now we take $\G$ to be an additive subgroup of $\Bbb{F}^m$ such that
$$\G\supset \{(j_1,...,j_{n_1},0,...,0)\mid j_1,...,j_{n_1}\in\Bbb{Z}\}\eqno(2.17)$$
and define
$$\vf_p=\zeta_p|_{\G},\;\;\vf_q\equiv 0\qquad\for\;\;p\in\ol{1,m},\;q\in\ol{m+1,n}.\eqno(2.18)$$
Then $\vf_1,...,\vf_n$ are additive group homomorphisms satisfying (2.6) and (2.7). For instance, we can take
$$\G=\sum_{j=1}^k((j/k,j/k,...,j/k)+\Bbb{Z}^m)\eqno(2.19)$$
for any positive integer $k$. When $n_1=0,\;m=n$ and $\G=\Bbb{Z}^n$,
$${\cal A}=\Bbb{F}[t_i,t_{n+i},t_i^{-1}\mid i\in\ol{1,n}]\eqno(2.20)$$
(cf. (2.10) and (2.11)) and 
$$\ptl_i=t_i\ptl_{t_i}+\ptl_{t_{n+i}}\qquad\for\;\;i\in\ol{1,n}.\eqno(2.21)$$

\section{Algebras of Type S}

In this section, we shall construct and prove a  new class of generalized simple Lie algebras of Special type. 

We shall use the same notations as in last section. For convenience, we denote
$$x_1^{\al}=x^{\al,\vec{0}},\;\;x_2^{\vec{i}}=x^{0,\vec{i}}\qquad\;\;\for\;\;\al\in\G,\;\vec{i}\in\vec{\cal J}.\eqno(3.1)$$ 
We shall restrict to a special case of $(\G,\vec{\vf})$. 
Let $\{\Delta_p\mid p\in\ol{1,n}\}$ be $n$ additive subgroups of $\Bbb{F}$.
Set 
$$\G=\Delta_1\times\Delta_2\times\cdots\times\Delta_n,\eqno(3.2)$$
where the addition on $\vec{\Delta}$ is defined componentwise. We use $\vec{\al}=(\al_1,\al_2,...,\al_n)$ to denote an element in $\G$ with $\al_p\in\Delta_p$. Moreover, for $p\in\ol{1,n}$, we take $\vf_p$ as
$$\vf_p(\vec{\al})=\al_p\qquad\for\;\;\vec{\al}\in \G.\eqno(3.3)$$
Furthermore, we denote 
$$\al_{[p]}=(0,..,\stackrel{p}{\al},0,...,0)\qquad\for\;\;\al\in\Delta_p.\eqno(3.4)$$
Now Condition (2.6) becomes
$$\Delta_p\neq\{0\}\qquad\mbox{if}\;\;{\cal J}_p=\{0\}\;\;\for\;\;p\in\ol{1,n}\eqno(3.5)$$
and Condition (2.7) holds automatically. 
\psp

Recall the algebra $({\cal A},\cdot)$ defined by  (2.10), (2.11) and $\{\ptl_1,...,\ptl_n\}$ defined by (2.12) and (2.14). Moreover, we assume that $n\geq 2$. 

Let $\vec{\rho},\vec{\sgm}\in\G$ be two fixed elements. We set
$$D_p=x_1^{(\sgm_p)_{[p]}}\ptl_p\qquad\for\;\;p\in\ol{1,n}\eqno(3.6)$$
(cf. (3.1) and (3.4)). Then $\{D_1,...,D_n\}$  are linearly independent and mutually commutative derivations of ${\cal A}$ satisfying (2.2).  Moreover, we set
$$D_{p,q}(u)=x_1^{\vec{\rho}}(D_q(u)D_p-D_p(u)D_q)\eqno(3.7)$$
for $p,q\in\ol{1,n},\;u\in{\cal A}$. By (3.1) in [O2],
\begin{eqnarray*}& &[D_{p,q}(u),D_{r,s}(v)]\\&=&D_{p,s}(x_1^{\vec{\rho}}D_q(u)D_r(v))+D_{q,r}(x_1^{\vec{\rho}}D_p(u)D_s(v))\\&&-D_{p,r}(x_1^{\vec{\rho}}D_q(u)D_s(v))-D_{q,s}(x_1^{\vec{\rho}}D_p(u)D_r(v))\hspace{6.1cm}(3.8)\end{eqnarray*}
for $p,q,r,s\in\ol{1,n}$ and $u,v\in{\cal A}$. We define:
$${\cal S}=\mbox{span}\:\{D_{p,q}(u)\mid p,q\in\ol{1,n},\;u\in{\cal A}\}\subset {\cal W}.\eqno(3.9)$$
In the spirit of (3.8), ${\cal S}$ forms a Lie subalgebra of ${\cal W}$, which we call a {\it   generalized Lie algebra of Special type}. The following result was proved by Osborn [O2].
\psp

{\bf Proposition 3.1}. {\it The Lie algebra} ${\cal S}^{(1)}=[{\cal S},{\cal S}]$ {\it is a simple Lie algebra when} $\Delta_p=\{0\}$ {\it or} ${\cal J}_p=\{0\}$ {\it for each} $p\in\ol{1,n}$.
\psp

Our main theorem in this section is:
\psp

{\bf Theorem 3.2}. {\it The Lie algebra} ${\cal S}$ {\it is simple if} $\Delta_p\neq\{0\}$ {\it and} ${\cal J}_p\neq \{0\}$ {\it for some} $p\in \ol{1,n}$.
\psp

We shall prove the theorem by establishing several lemmas. For convenience, we can assume:
$$\Delta_1\neq \{0\},\;\;\;{\cal J}_1\neq\{0\}.\eqno(3.10)$$
Set
$$\iota_{p,q}=(\rho_p+\sgm_p)_{[p]}+(\rho_q+\sgm_q)_{[q]},\;\;\vec{\gm}_{p,q}=(\gm_p)_{[p]}+(\gm_q)_{[q]}\qquad\for\;\;\vec{\gm}\in\G\eqno(3.11)$$
(cf. (3.4)).

 For any $(\vec{\al},\vec{i})$  and $(\vec{\be},\vec{j})\in\G\times \vec{\cal J}$, (3.8) implies the following identity
\begin{eqnarray*}& &[D_{p,q}(x^{\vec{\al},\vec{i}}),D_{p,q}(x^{\vec{\be},\vec{j}})]\\&=&(\al_q\be_p-\al_p\be_q)D_{p,q}(x^{\vec{\al}+\vec{\be}+\vec{\rho}+\vec{\sgm}_{p,q},\vec{i}+\vec{j}})+(i_qj_p-i_pj_q)D_{p,q}(x^{\vec{\al}+\vec{\be}+\vec{\rho}+\vec{\sgm}_{p,q},\vec{i}+\vec{j}-1_{[p]}-1_{[q]}})\\& &+(\al_qj_p-\be_qi_p)D_{p,q}(x^{\vec{\al}+\vec{\be}+\vec{\rho}+\vec{\sgm}_{p,q},\vec{i}+\vec{j}-1_{[p]}})\\& &+(\be_pi_q-\al_pj_q)D_{p,q}(x^{\vec{\al}+\vec{\be}+\vec{\rho}+\vec{\sgm}_{p,q},\vec{i}+\vec{j}-1_{[q]}})\hspace{7cm}(3.12)\end{eqnarray*}
(cf. (2.9)). 

\psp

In the rest of this section, we let $I$ be a nonzero ideal of ${\cal S}$. Moreover, by reindexing if necessary,  we can assume
$$\Delta_p\neq \{0\},\;\;\Delta_q=\{0\}\qquad\for\;\;p,q\in\ol{1,n},\;p\leq m,\;q>m,\eqno(3.13)$$
where $m\in \ol{1,n}$. Obviously, $m\geq 1$ by our assumption (3.10). A nonzero element 
$$u=\sum_{p=1}^n\sum_{\vec{\al}\in\G,\vec{i}\in\vec{\cal J}}a_{p,\vec{\al},\vec{i}}x^{\vec{\al},\vec{i}}\ptl_p\in {\cal W}\eqno(3.14)$$
is called a {\it homogeneous element of degree} $k$ if
$$\sum_{p=1}^mi_p=k\;\;\mbox{whenever}\;\;a_{q,\vec{\al},\vec{i}}\neq 0\;\;\for\;\;q\in\ol{1,n},\;(\vec{\al},\vec{i})\in\G\times \vec{\cal J}.\eqno(3.15)$$
For a homogeneous element $u\in {\cal W}$, we denote its degree by $\wp(u)$. The degree of $0$ is defined as $-1$. Furthermore, for $k\in\Bbb{N}$,  we define:
$${\cal W}_k=\mbox{span}\:\{\mbox{homogeneous element}\;u\in {\cal W}\mid \wp(u)\leq k\}.\eqno(3.16)$$
The leading term $u_{ld}$ of a nonzero element $u\in {\cal W}$ is a nonzero homogeneous element defined by 
$$u=u_{ld}+u'\qquad\mbox{with}\;\;u'\in {\cal W}_{k-1},\;\wp(u_{ld})=k\;\;\mbox{for some}\;\;k\in\Bbb{N}. \eqno(3.17)$$

\psp

{\bf Lemma 3.3}. {\it The subspace} $I_{1,2}=I\bigcap ({\cal A}\ptl_1+{\cal A}\ptl_2)\neq\{0\}$.
\psp

{\it Proof}. For any $u\in {\cal W}$, we write
$$u=u'+u^{\ast}\qquad\mbox{with}\;\;u'=\sum_{(\vec{\al},\vec{i})\in\G\times\vec{j}}x^{\vec{\al},\vec{i}}(a_{1,\vec{\al},\vec{i}}\ptl_1+a_{2,\vec{\al},\vec{i}}\ptl_2),\;\;u^{\ast}\in\sum_{p=3}^n{\cal A}\ptl_p,\eqno(3.18)$$
where $a_{1,\vec{\al},\vec{i}},a_{2,\vec{\al},\vec{i}}\in\Bbb{F}$. Moreover, we write the leading term:
$$(u^{\ast})_{ld}=\sum_{p=3}^n\sum_{(\vec{\al},\vec{i})\in\G\times\vec{\cal J}}b_{p,\vec{\al},\vec{i}}x^{\vec{\al},\vec{i}}\ptl_p\eqno(3.19)$$
with all $b_{p,\vec{\al},\vec{i}}\in \Bbb{F}$ and define
$$\imath(u)=|\{(p,\vec{\al},\vec{i})\mid b_{p,\vec{\al},\vec{i}}\neq 0,\;2< p\in\ol{1,n},\;\vec{\al}\in\G,\;\vec{i}\in\vec{\cal J}\}|.\eqno(3.20)$$
Set
$$\hat{k}=\mbox{min}\;\{\wp((u^{\ast})_{ld})\mid 0\neq u\in I\},\eqno(3.21)$$
where $\wp((u^{\ast})_{ld})$ is the degree of $(u^{\ast})_{ld}$ defined by (3.15). Our statement in the lemma is equivalent to $\hat{k}=-1$. Assume that $\hat{k}\geq 0$. We define
$$\imath=\min\{\imath(u)\mid 0\neq u\in I, \;\wp((u^{\ast})_{ld})=\hat{k}\}.\eqno(3.22)$$
If $\imath=0$, then we get a contradiction to (3.20). Now we assume $\imath>0$. Pick any $u\in I$ such that $\wp((u^{\ast})_{ld})=\hat{k}$ and  $\imath(u)=\imath$. We write $u$ as in (3.18) and (3.19).
\psp

{\it Case 1}. $\Delta_2\neq\{0\}$.
\psp

{\it Subcase 1}. There exist $b_{p,\vec{\al},\vec{i}}\neq 0$ and $b_{q,\vec{\be},\vec{j}}\neq 0$ such that $\al_1=0$ and  $\be_1\neq 0$.

 Pick any $0\neq \gm\in\Delta_2$. Note
$$D_{1,2}(x_1^{\gm_{[2]}})=\gm x_1^{\vec{\rho}+\vec{\sgm}_{1,2}+\gm_{[2]}}\ptl_1.\eqno(3.23)$$
Then 
$$v=[D_{1,2}(x_1^{\gm_{[2]}}),u]\in I,\;\;\wp((v^{\ast})_{ld})=\hat{k},\;\;0<\imath(v)<\imath(u)=\imath,\eqno(3.24)$$
which contradicts  (3.20). 
\psp

{\it Subcase 2}. All $\al_1=0$ for $b_{p,\vec{\al},\vec{i}}\neq 0$. 

Suppose that $b_{q,\vec{\be},\vec{j}}\neq 0$. We shall use (3.18) and (3.19). First we assume that $\Delta_q\neq \{0\}$. For any $0\neq \gm\in\Delta_q$, we have:
$$\gm^{-1}[u^{\ast},D_{1,q}(x_1^{\gm_{[q]}})]\equiv \sum_{\vec{\al}\in\G,\vec{i}\in\vec{\cal J}}(\gm+\rho_q+\sgm_q)b_{q,\vec{\al},\vec{i}}x^{\vec{\al}+\gm+\vec{\rho}+\vec{\sgm}_{1,q},\vec{i}}\ptl_1\;\;(\mbox{mod}\;{\cal W}_{\hat{k}-1}).\eqno(3.25)$$
Since the coefficients of $x^{\vec{\al},\vec{i}}\ptl_1$ and $x^{\vec{\al},\vec{i}}\ptl_2$ in $\gm^{-1}[u',D_{1,q}(x_1^{\gm_{[q]}})]$ (cf. (3.17)) are independent of $\gm$, only finite number of them are nonzero and $|\Delta_q|=\infty$, (3.25) enables us to choose $0\neq \gm\in \Delta_q$ such that
$$0\neq v=\gm^{-1}[u,D_{1,q}(x_1^{\gm_{[q]}})]\in I,\;\;\wp((v^{\ast})_{ld})<\hat{k},\eqno(3.26)$$
which contradicts (3.21).

Next we assume $\Delta_q= \{0\}$. Then ${\cal J}_q=\Bbb{N}$ by our earlier assumption. 
For any $0\neq l\in\Bbb{N}$, we have:
$$l^{-1}[u^{\ast},D_{1,q}(x_2^{l_{[q]}})]\equiv \sum_{\vec{\al}\in\G,\vec{i}\in\vec{\cal J}}(l-1)b_{q,\vec{\al},\vec{i}}x^{\vec{\al}+\vec{\rho}+\vec{\sgm}_{1,q},\vec{i}+(l-2)_{[q]}}\ptl_1\;\;(\mbox{mod}\;{\cal W}_{\hat{k}-1}).\eqno(3.27)$$
Since the coefficients of $x^{\vec{\al},\vec{i}}\ptl_1$ and $x^{\vec{\al},\vec{i}}\ptl_2$ in $l^{-1}[u',D_{1,q}(x_2^{l_{[q]}})]$ are independent of $l$ and only finite number of them are nonzero, (3.27) enables us to choose $0\neq l\in \Bbb{N}$ such that
$$0\neq v=l^{-1}[u,D_{1,q}(x_2^{l_{[q]}})]\in I,\;\;\wp((v^{\ast})_{ld})<\hat{k},\eqno(3.28)$$
which contradicts (3.21).
\psp

We have the same conclusion if $\al_2=0$ for some $b_{p,\vec{\al},\vec{i}}\neq 0$.
\psp

{\it Subcase 3}. All $\al_1\neq 0$ and $\al_2\neq 0$, whenever $b_{p,\vec{\al},\vec{i}}\neq 0$.
\psp

Assume that $b_{q,\vec{\be},\vec{j}}\neq 0$. If $\be_1\neq -(\rho_1+\sgm_1)$, then
$$v=[D_{1,2}(x_1^{-(\be_1+\rho_1+\sgm_1)_{[1]}}),u]\in I,\;\;\wp((v^{\ast})_{ld})=\hat{k},\;\;\imath(v)=\imath(u)\eqno(3.29)$$
and $v$ has a term
$$(\be_1+\rho_1+\sgm_1)\be_2b_{q,\vec{\be},\vec{j}}x^{\vec{\be}-(\be_1)_{[1]}+\vec{\rho}-(\rho_1)_{[1]}+(\sgm_2)_{[2]},\vec{j}}\ptl_q\qquad\mbox{if}\;\;\Delta_2\neq\{0\}.\eqno(3.30)$$
Replacing $u$ by $v$, we go back to the above Subcases 1 and 2. If $\be_2+\rho_2+\sgm_2=0$, we choose any $0\neq \gm\in \Delta_1$. Replacing $u$ by $[D_{1,2}(x_1^{\gm_{[1]}}),u]$, we get the same situation as in the above Subcases 1 and 2 with indices 1 and 2 exchanged. Now we assume that $\be_2+\rho_2+\sgm_2\neq 0$ and $\be_1= -(\rho_1+\sgm_1)$. Pick any $\gm\in\Delta_1$ such that $\gm\neq 0,\;-(\rho_1+\sgm_1)$. We have:
$$v=[D_{1,2}(x_1^{-(\gm+\rho_1+\sgm_1)_{[1]}}),[D_{1,2}(x_1^{\gm_{[1]}}),u]]\in (I\bigcap {\cal W}_k)\setminus {\cal W}_{k-1},\;\;\imath(v)=\imath(u).\eqno(3.31)$$
Replacing $u$ by $v$, we again go back to the above Subcases 1 and 2. 

In summary, we always get a contradiction if $\Delta_2\neq\{0\}$.
\psp

{\it Case 2}. $\Delta_2=\{0\}$.
\psp

According to our assumption, $m=1$ in this case (cf. (3.13)). Thus ${\cal J}_2=\cdots={\cal J}_n=\Bbb{N}$. 

\pse

{\it Subcase 1}. There exist $b_{p,\vec{\al},\vec{i}}\neq 0$ and $b_{q,\vec{\be},\vec{j}}\neq 0$ such that $\al_1=0$ and $\be_1\neq 0$.  

Note
$$D_{1,2}(x_2^{1_{[2]}})=x_1^{\vec{\rho}+(\sgm_1)_{[1]}}\ptl_1.\eqno(3.32)$$
$$v=[D_{1,2}(x_2^{1_{[2]}}),u]\in I,\;\;\wp((v^{\ast})_{ld})=\hat{k},\;\;0<\imath(v)<\imath(u_{ld})=\imath,\eqno(3.33)$$
which contradicts (3.22). 
\psp

{\it Subcase 2}. All $\al_1=0$ for $b_{p,\vec{\al},\vec{i}}\neq 0$. 

This is the same as the situation that $\Delta_q= \{0\}$ in Subcase 2 of Case 1.
\psp

We have the same conclusion if $i_2=0$ for some $b_{p,\vec{\al},\vec{i}}\neq 0$.
\psp

{\it Subcase 3}. All $\al_1\neq 0,\;i_2\neq 0$ whenever $b_{p,\vec{\al},\vec{i}}\neq 0$.

Assume that $b_{q,\vec{\be},\vec{j}}\neq 0$. If $\be_1\neq -(\rho_1+\sgm_1)$, then we have (3.29), and $v$ has a term
$$(\be_1+\rho_1+\sgm_1)j_2b_{q,\vec{\be},\vec{j}}x^{\vec{\be}-(\be_1)_{[1]}+\vec{\rho}-(\rho_1)_{[1]},\vec{j}-1_{[2]}}\ptl_q.\eqno(3.34)$$
Replacing $u$ by $v$, we go back to the above Subcases 1 and 2. If $j_2=1$, we choose any $0\neq \gm\in \Delta_1$. Replacing $u$ by $[D_{1,2}(x_1^{\gm_{[1]}}),u]$, we get the same situation as in the above with indices 1 and 2 exchanged. Now we assume that $j_2>1$ and $\be_1= -(\rho_1+\sgm_1)$. Pick any $\gm\in\Delta_1$ such that $\gm\neq 0,\;-(\rho_1+\sgm_1)$. We have (3.31). 
Thus, we always get a contradiction if $\Delta_2=\{0\}$.
\psp

Therefore, we get a contradiction if $\hat{k}\geq 0$. This completes the proof of the lemma.$\qquad\Box$
\psp

For $\vec{i}\in \vec{\cal J}$, we define
$$|\vec{i}|_{1,2}=\left\{\begin{array}{ll}i_1&\mbox{if}\;\;\Delta_2=\{0\},\\i_1+i_2&\mbox{if}\;\;\Delta_2\neq\{0\}.\end{array}\right.\eqno(3.35)$$
Moreover, for 
$$u=\sum_{(\al,\vec{i})\in\G\times \vec{\cal J}}x^{\al,\vec{i}}(a_{1,\al,\vec{i}}\ptl_1+a_{2,\al,\vec{i}}\ptl_2)\in{\cal W},\eqno(3.36)$$
we let
$$\wp_{1,2}(u)=\max\:\{|\vec{i}|_{1,2}\mid a_{1,\al,|\vec{i}|}\neq 0\;\mbox{or}\;a_{2,\al,|\vec{i}|}\neq 0\}.\eqno(3.37)$$
In addition, we put $\wp_{1,2}(0)=-1$. Set
$${\cal W}_{1,2}^k=\{u\in {\cal A}\ptl_1+{\cal A}\ptl_2\mid \wp_{1,2}(u)\leq k\}\qquad\for\;\;k\in\Bbb{N}.\eqno(3.38)$$
\psp

{\bf Lemma 3.4}. {\it If} $(I_{1,2}\bigcap {\cal W}_{1,2}^k)\setminus {\cal W}_{1,2}^{k-1}\neq \emptyset$ {\it for some} $k\in \Bbb{N}$, {\it then there exists an element} $u\in (I_{1,2}\bigcap {\cal W}_{1,2}^k)\setminus {\cal W}_{1,2}^{k-1}$ {\it such that}
$$u\equiv\sum_{\vec{\al}\in\G,\vec{i}\in\vec{\cal J}}c_{\vec{\al},\vec{i}}D_{1,2}(x^{\vec{\al},\vec{i}})\;\;(\mbox{\it mod}\;{\cal W}_{1,2}^{k-1})\eqno(3.39)$$
{\it with}
$$(\al_1,\al_2)=(\kappa_1,\kappa_2)\;\;\mbox{\it whenever}\;\;c_{\vec{\al},\vec{i}}\neq 0\eqno(3.40)$$
{\it if} $\Delta_2\neq 0$, {\it or} 
$$(\al_1,i_2)=(\kappa_1,\ves)\;\;\mbox{\it whenever}\;\;c_{\vec{\al},\vec{i}}\neq 0\eqno(3.41)$$
{\it if} $\Delta_2=\{0\}$, {\it where} $\kappa_p\in\Delta_p$ {\it and} $\ves\in\Bbb{N}$.
\psp

{\it Proof}. Assume $(I_{1,2}\bigcap {\cal W}_{1,2}^k)\setminus {\cal W}_{1,2}^{k-1}\neq \{0\}$. For any $u\in (I_{1,2}\bigcap {\cal W}_{1,2}^k)\setminus {\cal W}_{1,2}^{k-1}$, we write
$$u\equiv\sum_{\vec{\al}\in\G,\vec{i}\in\vec{\cal J}}a_{\vec{\al},\vec{i}}D_{1,2}(x^{\vec{\al},\vec{i}})\;\;(\mbox{mod}\;{\cal W}_{1,2}^{k-1})\eqno(3.42)$$
with
$$D_{1,2}(x^{\vec{\al},\vec{i}})\in {\cal W}_{1,2}^k\setminus {\cal W}_{1,2}^{k-1}\qquad\mbox{whenever}\;\;a_{\vec{\al},\vec{i}}\neq 0.\eqno(3.43)$$
Set
$$\ell(u)=\left\{\begin{array}{ll}|\{(\al_1,\al_2)\mid a_{\vec{\al},\vec{i}}\neq 0\}|&\mbox{if}\;\;\Delta_2\neq \{0\},\\ |\{(\al_1,i_2)\mid a_{\vec{\al},\vec{i}}\neq 0\}|&\mbox{if}\;\;\Delta_2= \{0\}.\end{array}\right.\eqno(3.44)$$
Moreover, we let
$$\ell=\min\{\ell(u)\mid u\in (I_{1,2}\bigcap {\cal W}_{1,2}^k)\setminus {\cal W}_{1,2}^{k-1}\}.\eqno(3.45)$$
We want to prove that $\ell=1$.

Assume $\ell>1$. We pick $u\in (I_{1,2}\bigcap {\cal W}_{1,2}^k)\setminus {\cal W}_{1,2}^{k-1}$ such that $\ell(u)=\ell$. Write $u_{ld}$ by (3.42) and (3.43). 
\psp

{\it Case 1}. $\Delta_2\neq \{0\}$.
\psp

{\it Subcase 1}. There exist $a_{\vec{\al},\vec{i}}\neq 0$ and $a_{\vec{\be},\vec{j}}\neq 0$ such that
$(\be_1,\be_2)\not \in \Bbb{F}(\al_1,\al_2)$.

In this situation, we have
$$v=[D_{1,2}(x^{\vec{\al}}_1),u]\in (I_{1,2}\bigcap {\cal W}_{1,2}^k)\setminus {\cal W}_{1,2}^{k-1},\;\;\;\ell(v)<\ell(u)\eqno(3.46)$$
by (3.11), which contradicts (3.45).
\psp

{\it Subcase 2}. Assume $a_{\vec{\al},\vec{i}}\neq 0$. If $a_{\vec{\be},\vec{j}}\neq 0$, then $(\be_1,\be_2)\in\Bbb{F}(\al_1,\al_2)$. 

Since $\ell>1$, there exists $b_{\vec{\be},\vec{j}}\neq 0$ such that $(\be_1,
\be_2)\neq (\al_1,\al_2)$. Pick any $(\gm_1,\gm_2)\in\Delta_1\times \Delta_2$ such that $(\gm_1,\gm_2),(\gm_1+\rho_1+\sgm_1,\gm_2+\rho_2+\sgm_2)\not\in \Bbb{F}(\al_1,\al_2)$. Then 
$$v=[D_{1,2}(x_1^{(\gm_1)_{[1]}+(\gm_2)_{[2]}}),u]\in (I_{1,2}\bigcap {\cal W}_{1,2}^k)\setminus {\cal W}_{1,2}^{k-1},\;\;\;\ell(v)=\ell(u).\eqno(3.47)$$
Replacing $u$ by this $v$, we go back to the above Subcase 1.
\psp

{\it Case 2}. $\Delta_2=\{0\}$. In this case, ${\cal J}_2=\Bbb{N}$.
\psp

{\it Subcase 1}. There exist $a_{\vec{\al},\vec{i}}\neq 0$ and $a_{\vec{\be},\vec{j}}\neq 0$ such that
$(\be_1,j_2)\not \in \Bbb{F}(\al_1,i_2)$.

In this situation, we have (3.46) with $x_1^{\vec{\al}}$ replaced by $x^{\vec{\al},(i_2)_{[2]}}$ , which contradicts (3.45).

{\it Subcase 2}. Assume $a_{\vec{\al},\vec{i}}\neq 0$. If $a_{\vec{\be},\vec{j}}\neq 0$, then $(\be_1,j_2)\in\Bbb{F}(\al_1,i_2)$. 

Since $\ell>1$, there exists $a_{\vec{\be},\vec{j}}\neq 0$ such that $(\be_1,
j_2)\neq (\al_1,i_2)$. Pick any $(\gm,l)\in\Delta_1\times \Bbb{N}$ such that $(\gm,l),(\gm+\rho_1+\sgm_1,l)\not\in \Bbb{F}(\al_1,i_2)$. Then 
$$v=[D_{1,2}(x^{\gm_{[1]},l_{[2]}}),u]\in (I_{1,2}\bigcap {\cal W}_{1,2}^k)\setminus {\cal W}_{1,2}^{k-1},\;\;\;\ell(v)=\ell(u)\eqno(3.48)$$
by (3.11). Replacing $u$ by this $v$, we go back to the above Subcase 1.
\psp

Therefore, $\ell=1$, which implies our conclusion in the lemma.$\qquad\Box$
\psp

Set
$$\td{k}=\min\{k\mid I_{1,2}\bigcap {\cal W}_{1,2}^k\neq \{0\}\},\eqno(3.49)$$
which is well-defined by Lemma 3.3.

\psp

{\bf Lemma 3.5}. {\it There exists an element} $u\in I_{1,2}$ such that 
$$u=\sum_{\vec{\al}\in\G,\vec{i}\in\vec{\cal J}}d_{\vec{\al},\vec{i}}D_{1,2}(x^{\vec{\al},\vec{i}})\eqno(3.50)$$
{\it with} 
$$\al_1=0,\;\;i_1=0\;\;\;\;(\al_2,i_2)=(\kappa_2,\ves)\qquad\mbox{\it whenever}\;\;d_{\vec{\al},\vec{i}}\neq 0,\eqno(3.51)$$
{\it for any} $\kappa_2\in\Delta_2$ {\it and} $\ves\in\Bbb{N}$ {\it such that} 
$$\ves =0,\;\;\kappa_2\neq 0,\;\rho_2+\sgm_2\qquad\mbox{\it if}\;\Delta_2\neq \{0\};\qquad\;\;\ves\neq 0\qquad\mbox{if}\;\Delta_2=\{0\}.\eqno(3.52)$$
\psp

{\it Proof}. Take $u$ to be an element as in the statement of  Lemma 3.4 with $k=\td{k}$.
\psp

{\it Case 1}. $\Delta_2\neq \{0\}$.
\psp

  Assume that $\td{k}>0$. If $(\kappa_1,\kappa_2)=(0,0)$, then $0\leq \wp_{1,2}([D_{1,2}(x_1^{(\gm_1)_{[1]}+(\gm_2)_{[2]}}),u])<\td{k}$ for any $0\neq \gm_1\in\Delta_1,\;\;0\neq \gm_2\in\Delta_2$, which contradicts (3.49). Assume $(\kappa_1,\kappa_2)\neq(0,0)$.  Note that
\begin{eqnarray*}& &[D_{1,2}(x_1^{(\kappa_1)_{[1]}+(\kappa_2)_{[2]}}),u]\\&\equiv& \sum_{\vec{\al}\in\G,\vec{i}\in\vec{\cal J}}c_{\vec{\al},\vec{i}}[\kappa_2 i_1 D_{1,2}(x^{\vec{\al}+\vec{\rho}+\vec{\sgm}_{1,2}+(\kappa_1)_{[1]}+(\kappa_2)_{[2]},\vec{i}-1_{[1]}})\\& &-\kappa_1i_2D_{1,2}(x^{\vec{\al}+\vec{\rho}+\vec{\sgm}_{1,2}+(\kappa_1)_{[1]}+(\kappa_2)_{[2]},\vec{i}-1_{[2]}})]\\& &+
\sum_{(\vec{\be},\vec{j})\in\G\times\vec{\cal J};\be_p\neq 2\kappa_p+\sgm_p+\rho_p,\:p=1\;\mbox{\tiny or}\;2}a'_{\vec{\be},\vec{j}}D_{1,2}(x^{\vec{\be},\vec{j}})\;\;(\mbox{mod}\;{\cal W}_{1,2}^{\td{k}-2}).\hspace{3.5cm}(3.53)\end{eqnarray*}
Considering $[D_{1,2}(x_1^{\vec{\gm}}),u]\in I$ if necessary, we can take any $(\kappa_1,\kappa_2)\neq (\rho_1+\sgm_1,\rho_2+\sgm_2)$  without changing $\vec{i}$ for all $c_{\vec{\al},\vec{i}}\neq 0$ in (3.39) by (3.11) and the proof of Proposition 3.6 in [O2] (also cf. (3.60)-(3.67)). Choose $u$ such that $\kappa_1\neq 0,\;\kappa_2=0$ if $i_2\neq 0$ and $\kappa_2\neq 0,\;\kappa_1=0$ otherwise. For such $u$,  $[D_{1,2}(x_1^{(\kappa_1)_{[1]}+(\kappa_2)_{[2]}}),u]\in (I\bigcap {\cal W}_{1,2}^{\td{k}-1})\setminus {\cal W}_{1,2}^{\td{k}-2}$, which contradicts (3.49). So $\td{k}=0$. Considering $[D_{1,2}(x_1^{\vec{\gm}}),u]\in I$  we can obtain $\kappa_1=0$ and $\kappa_2\neq 0,\rho_2+\sgm_2$  by (3.11) and the proof of Proposition 3.6 in [O2] (also cf. (3.60)-(3.67)).
\psp

{\it Case 2}. $\Delta_2=\{0\}$.
\psp

In this case, ${\cal J}_2=\Bbb{N}$ by (2.6). Assume that $\td{k}>0$. If $(\kappa_1,\ves)=(0,0)$, then we have $0\leq \wp_{1,2}([D_{1,2}(x_1^{(\gm_1)_{[1]},1_{[2]}}),u])<\td{k}$ for any $0\neq \gm_1\in\Delta_1,\;\;0\neq l\in\Bbb{N}$, which contradicts (3.49). Assume $(\kappa_1,\ves)\neq(0,0)$. Considering $[D_{1,2}(x^{\gm_{[1]},l_{[2]}}),u]\in I$ if necessary, we can take any $(\kappa_1,\ves)\neq (\rho_1+\sgm_1,0)$  without changing $\vec{i}$ for all $c_{\vec{\al},\vec{i}}\neq 0$ in (3.39) by (3.11) and the proof of Proposition 3.6 in [O2] (also cf. (3.60)-(3.67)). In particular, we can take $\kappa_1=0$ and $\ves\neq0$. Assume that $i_1\neq 0$ for some $c_{\vec{\al},\vec{i}}\neq 0$ in (3.39). Since
\begin{eqnarray*}& &[D_{1,2}(x_2^{1_{[2]}}),u]\\&\equiv& \sum_{\vec{\al}\in\G,\vec{i}\in\vec{\cal J}}c_{\vec{\al},\vec{i}}i_1 D_{1,2}(x^{\vec{\al}+\vec{\rho}+\vec{\sgm}_{1,2},\vec{i}-1_{[1]}})\\& &+
\sum_{(\vec{\be},\vec{j})\in\G\times \vec{\cal J};\be_1\neq \rho_1+\sgm_1\;\mbox{\tiny or}\;j_2\neq\ves}b_{\vec{\be},\vec{j}}D_{1,2}(x^{\vec{\be},\vec{j}})\;\;(\mbox{mod}\;{\cal W}_{1,2}^{\td{k}-2}),\hspace{4.6cm}(3.54)\end{eqnarray*}
we have $[D_{1,2}(x_2^{1_{[2]}}),u]\in (I\bigcap {\cal W}_{1,2}^{\td{k}-1})\setminus {\cal W}_{\td{k}-2}$, which contradicts (3.49). Thus $i_1=0.\qquad\Box$
\psp

{\bf Lemma 3.6}. {\it We have} $D_{1,n}(x_1^{(\al_1)_{[1]}})\in I$ {\it for some} $0\neq \al_1\in\Delta_1$.
\psp

{\it Proof}. For any $p\in\ol{1,n}$, we define:
$$\Psi_p=\left\{\begin{array}{ll}(\Delta_p,0)&\mbox{if}\;\;\Delta_p\neq\{0\},\\(0,\Bbb{N})&\mbox{if}\;\;\Delta_p=\{0\}.\end{array}\right.\eqno(3.55)$$
For $p>3$, we define $I_{1,p}$ to be the subset of $I$ of the elements of the form:
$$u=\sum_{\vec{\al}\in\G,\vec{i}\in\vec{\cal J}}a_{\vec{\al},\vec{i}}D_{1,p}(x^{\vec{\al},\vec{i}})\eqno(3.56)$$
with $a_{\vec{\al},\vec{i}}\in\Bbb{F}$ and
$$(\al_q,i_q)=(\kappa_q,\ves_q)\qquad \for\;\;2\leq q\leq p-1\qquad\mbox{whenever}\;\;a_{\vec{\al},\vec{i}}\neq 0,\eqno(3.57)$$
where
$$(\kappa_q,\ves_q)\in\Psi_q\qquad \for\;\;2\leq q\leq p-1\eqno(3.58)$$
are fixed for each $u$. 

Let $u$ be an element as Lemma 3.5 (cf. (3.50)). By considering $[u,D_{1,2}(x^{\vec{\al},\vec{i}})]$ with $i_1=0,\;(\al_2,i_2)\in\Psi_2$ (cf. (3.11)) and the proof of Proposition 3.6 in [O2]), we can assume that $u$ satisfies (3.50)-(3.52) and 
$$(\al_3,i_3)\neq (0,0)\qquad\mbox{whenever}\;\;d_{\vec{\al},\vec{i}}\neq 0.\eqno(3.55)$$
Note that such an element $u\in I_{1,3}$. Thus $I_{1,3}\neq\{0\}$. Replacing the index 2 by 3 in the proofs of Lemmas 3.4 and 3.5, we can prove that $I_{1,4}\neq\{0\}$. Continuing this process, we can prove that $I_{1,n}\neq\{0\}$ by induction. Replacing the index 2 by $n$ in the proofs of Lemmas 3.4 and 3.5, we get $u=D_{1,n}(x^{\vec{\al},\vec{i}})\in I$ with $\al_1=0,\;i_1=0,\;(\al_p,i_p)\in\Psi_p$ for $2\leq p\in\ol{1,n}$ and $(\al_n,i_n)\neq (0,0)$. Repeating considering $[u,D_{1,n}(x^{\vec{\be},\vec{j}})]$ with $j_1=0,\;(\be_p,j_p)\in\Psi_p$ for $2\leq p\in\ol{1,n}$ by (3.11) and the proof of Proposition 3.6 in [O2] (also cf. (3.60)-(3.67)) , we can obtain $D_{1,n}(x_1^{(\al_1)_{[1]}})\in I$ for some $0\neq \al_1\in\Delta_1$. $\qquad\Box$
\psp

{\it Proof of Theorem 3.2}.
\psp

 We want to prove that $I={\cal S}$. By Lemma 3.6 and reindexing, we can assume that $D_{1,2}(x_1^{(\al_1)_{[1]}})\in I$ for convenience.
\psp

{\it Case 1}. $\Delta_2\neq 0$.
\psp

For any $\vec{\be}\in \G$ such that $\be_2\neq 0$, we have:
$$I\ni [D_{1,2}(x_1^{\al_{[1]}}),D_{1,2}(x_1^{\vec{\be}})]=-\al\be_2D_{1,2}(x_1^{\al_{[1]}+\vec{\be}+\vec{\rho}+\vec{\sgm}_{1,2}})\eqno(3.60)$$
by (3.11). So we obtain
$$D_{1,2}(x_1^{\al_{[1]}+\vec{\be}+\vec{\rho}+\vec{\sgm}_{1,2}})\in I.\eqno(3.61)$$
This implies that
$$D_{1,2}(x_1^{\vec{\be}+\iota_{1,2}})\in I\qquad\for\;\;\vec{\be}\in\G,\;\be_2\neq 0.\eqno(3.62)$$
If $\rho_2+\sgm_2= 0$, then $\vec{\gm}_{1,2}-2\iota_{1,2}-\vec{\be}_{1,2}$ and $\vec{\be}_{1,2}+\iota_{1,2}$ are linearly independent for any $\vec{\be},\vec{\gm}\in\G$ such that $\gm_1\neq \rho_1+\sgm_1,\;\gm_2=0$ and $\be_2\neq 0$. For such  $\vec{\be}$ and $\vec{\gm}$ , since
$$[D_{1,2}(x_1^{\vec{\gm}-\vec{\be}-\vec{\rho}-\vec{\sgm}_{1,2}-\iota_{1,2}}),D_{1,2}(x_1^{\vec{\be}+\iota_{1,2}})]\in I,\eqno(3.54)$$
we have
$$ D_{1,2}(x_1^{\vec{\gm}})\in I\eqno(3.64)$$
by (3.11). If $\rho_2+\sgm_2\neq 0$, then $ D_{1,2}(x_1^{2\iota_{1,2}})\in I$ by (3.63). Then for any $\gm\in\G$ such that $\gm_1\neq 0$ and $\gm_2=0$, since
$$[D_{1,2}(x_1^{\vec{\gm}-2\iota_{1,2}}), D_{1,2}(x_1^{2\iota_{1,2}})]=2\gm_1(\rho_2+\sgm_2)D_{1,2}(x_1^{\vec{\gm}+\vec{\rho}+\vec{\sgm}_{1,2}})\in I,\eqno(3.65)$$
we have
$$D_{1,2}(x_1^{\vec{\be}+\iota_{1,2}})\in I\qquad\for\;\;\vec{\be}\in\G,\;\be_1\neq 0.\eqno(3.66)$$
Combining (3.63) and (3.66), we get
$$D_{1,2}(x_1^{\vec{\be}})\in I\qquad\for\;\;\vec{\be}\in\G,\;\vec{\be}_{1,2}\neq\iota_{1,2}.\eqno(3.67)$$
 For any $\vec{\be}\in\G$ such that $\be_1,\be_2\neq 0$ and $\vec{\be}_{1,2}\neq \iota_{1,2}$,
$$[D_{1,2}(x^{\vec{\be},1_{[1]}}),D_{1,2}(x_1^{-\vec{\be}_{1,2}})]=\be_2D_{1,2}(x^{\vec{\be}-\vec{\be}_{1,2}+\vec{\rho}+\vec{\sgm}_{1,2}})\in I\eqno(3.68)$$
by (3.11). Therefore,
$$D_{1,2}(x_1^{\vec{\be}})\in I\qquad\mbox{for any}\;\;\vec{\be}\in\G.\eqno(3.69)$$
\psp

{\it Case 2}. $\Delta_2=0$.
\psp

According to our assumption, ${\cal J}_2=\Bbb{N}$ in this case. Since  for any $\vec{\be}\in \G$ and $0<j\in\Bbb{N}$, we have
$$I\ni [D_{1,2}(x_1^{\al_{[1]}}),D_{1,2}(x^{\vec{\be},j_{[2]}})]=-\al jD_{1,2}(x^{\al_{[1]}+\vec{\be}+\vec{\rho}+\vec{\sgm}_{1,2},(j-1)_{[2]}}).\eqno(3.70)$$
Thus
$$D_{1,2}(x^{\vec{\be},j_{[2]}})\in I\qquad\mbox{for any}\;\;\vec{\be}\in\G,\;j\in\Bbb{N}.\eqno(3.71)$$
\psp

We now want to prove $I={\cal S}$ by (3.69) and (3.71). Let $2<r\in\ol{1,n}$. If $\Delta_2\neq\{0\}$, we have:
$$[D_{1,2}(x_1^{\be_{[2]}}),D_{2,r}(x_1^{\vec{\gm}})]=\be\gm_2D_{1,r}(x_1^{\be_{[2]}+\vec{\gm}+\vec{\rho}+(2\sgm_2)_{[2]}})-\be\gm_rD_{1,2}(x_1^{\be_{[2]}+\vec{\gm}+\vec{\rho}+\vec{\sgm}_{2,r}})\in I\eqno(3.72)$$
for any $\be\in\Delta_2$ and $\vec{\gm}\in\G$ by (3.7). Thus by (3.69),
$$D_{1,r}(x_1^{\vec{\al}})\in I\qquad\mbox{for any}\;\;\vec{\al}\in \G\eqno(3.73)$$
if $\Delta_2\neq\{0\}$. If $\Delta_2=\{0\}$, then ${\cal J}_2=\Bbb{N}$ and (3.71) holds. In this case, we have:
$$[D_{1,2}(x_2^{1_{[2]}}), D_{2,r}(x^{\vec{\be},1_{[2]}})]=D_{1,r}(x_1^{\vec{\be}+\vec{\rho}})-\be_r D_{1,2}(x^{\vec{\be}+\vec{\rho}+(\sgm_r)_{[r]},1_{[2]}})\in I\eqno(3.70)$$
for any $\vec{\be}\in \G$. Thus (3.73) holds again.

Now we let $1<r\in\ol{1,n}$. If ${\cal J}_r\neq\{0\}$, the for any $\vec{\be}\in\G$ and $\vec{j}\in{\cal J}$, we have:
$$[D_{1,r}(x_2^{1_{[r]}}),D_{1,r}(x^{\vec{\be},\vec{j}})]=\be_1 D_{1,r}(x^{\vec{\be}+\vec{\rho}+\vec{\sgm}_{1,r},\vec{j}})+j_1D_{1,r}(x^{\vec{\be}+\vec{\rho}+\vec{\sgm}_{1,r},\vec{j}-1_{[1]}})\in I\eqno(3.75)$$
by (3.11) and (3.71). Moreover, by (3.75) and induction on $j_1$, we can prove that
$$D_{1,r}(x^{\vec{\be},\vec{j}})\in I\qquad\mbox{for any}\;\;\vec{\be}\in\G,\;\vec{j}\in\vec{\cal J}.\eqno(3.76)$$
If ${\cal J}_r=\{0\}$, then $\Delta_r\neq \{0\}$ by our assumption. Exchanging the positions of $1$ and $r$, we can prove (3.76) because $D_{1,r}(x^{\al_{[r]}})\in I$ for $\al\in \Delta_r$ by (3.65) and (3.73).

Now for any $1<r,s\in\ol{1,n},\;\vec{\be},\vec{\gm}\in\G$ and $\vec{j}\in\vec{\cal J}$, we have:
\begin{eqnarray*}[D_{1,r}(x^{\vec{\be},\vec{j}}),D_{1,s}(x_1^{\vec{\gm}})]&\equiv &\be_1\gm_1D_{s,r}(x^{\vec{\be}+\vec{\gm}+\vec{\rho}+(2\sgm_1)_{[1]},\vec{j}})\\& &+j_1\gm_1D_{s,r}(x^{\vec{\be}+\vec{\gm}+\vec{\rho}+(2\sgm_1)_{[1]},\vec{j}-1_{[1]}})\;\;(\mbox{mod}\:I).\hspace{3cm}(3.77)\end{eqnarray*}
by (3.8). By (3.77) and induction on $j_1$, we can prove that
$$D_{s,r}(x^{\vec{\be},\vec{j}})\in I\qquad\mbox{for any}\;\;\vec{\be}\in\G,\;\vec{j}\in\vec{\cal J}.\eqno(3.78)$$
Therefore, $I={\cal S}.\qquad\Box$

\section{Algebras of Type H}

In this section, we shall construct and prove a new class of generalized  simple Lie algebras of Hamiltonian type. 

All the notations and assumptions except (2.6) and (2.7) are the same as in Section 2. Moreover, we allow $n=0$.
Assume 
$$n=2m\qquad\mbox{for some}\;\;m\in\Bbb{N}\eqno(4.1)$$
and allow $m=0$. Moreover, we shall use the settings in (2.5), (2.8)-(2.14). 
Let $m_1\in\Bbb{N}$ such that $m_1\leq m$.  We view $\G$ as a $\Bbb{Z}$-module. Let $\phi(\cdot,\cdot):\G\times \G\rightarrow \Bbb{F}$ be a skew-symmetric $\Bbb{Z}$-bilinear form. We shall replace the assumptions (2.6) and (2.7) as follows. First, we assume
$$\vf_p(\G)+{\cal J}_p\neq\{0\}\qquad\for\;\;p\in\ol{1,n};\eqno(4.2)$$
$$\vf_p\neq 0\;\;\mbox{or}\;\;\vf_{m+p}\neq 0\qquad\for\;\;p\in\ol{1,m_1};\eqno(4.3)$$
$$\vf_{m+q}=0\;\mbox{if}\;{\cal J}_q=\{0\}\;\mbox{and}\;\vf_q=0\;\mbox{if}\;{\cal J}_{m+q}=\{0\}\eqno(4.4)$$
for $q\in\ol{m_1+1,m}$. Set
$$\mho=\{p,m+p\mid p\in\ol{1,m_1},\;\vf_p\neq 0,\;\vf_{m+p}\neq 0\}.\eqno(4.5)$$
Furthermore, we assume
$$(\bigcap_{p\neq q\in\ol{1,n}}\kn_{\vf_q})\setminus\kn_{\vf_p}\neq\emptyset\qquad\mbox{if}\;\;\vf_p\neq 0,\;\;p\in\ol{m_1+1,m}\bigcup \ol{m+m_1+1,n};\eqno(4.6)$$
$$\rad_{\phi}\bigcap (\bigcap_{p\neq q\in\ol{1,n}}\kn_{\vf_q})\setminus\kn_{\vf_p}\neq\emptyset\qquad\mbox{if}\;\;\vf_p\neq 0,\;\;p\in\ol{1,m_1}\bigcup \ol{m+1,m+m_1};\eqno(4.7)$$
$$\{\al\in\G\mid \phi(\al,\be)=0\;\for\;\be\in\bigcap_{q\in\mho}\kn_{\vf_q}\}\bigcap(\bigcap_{p=1}^{n}\kn_{\vf_p})=\{0\}.\eqno(4.8)$$
 We choose fixed elements
$$0\neq\sgm_p=\sgm_{m+p}\in\rad_{\phi}\bigcap(\bigcap_{p,m+p\neq q\in\ol{1,n}}\kn_{\vf_q})\qquad\for\;\;p\in\ol{1,m_1}.\eqno(4.9)$$
Set
$$\sgm=\sum_{p=1}^{m_1}\sgm_p.\eqno(4.10)$$

For any $\al\in \G$, we set
$${\cal A}_{\al}=\mbox{span}\:\{x^{\al,\vec{j}}\mid \vec{j}\in \vec{\cal J}\}.\eqno(4.11)$$
Moreover, we use notions in (3.1). 
Define an algebraic operation $[\cdot,\cdot]$ on ${\cal A}$ by:
\begin{eqnarray*}[u,v]&=&\sum_{q=m_1+1}^{m}[\ptl_q(u)\hat{\ptl}_{m+q}(v)+\hat{\ptl_q}(u)\ptl_{\vf_{m+q}}(v)-\hat{\ptl}_{m+q}(u)\ptl_q(v)-\ptl_{\vf_{m+q}}(u)\hat{\ptl}_q(v)]\\& &+\sum_{p=1}^{m_1}x_1^{\sgm_p}[\ptl_p(u)\ptl_{m+p}(v)-\ptl_{m+p}(u)\ptl_p(v)]+\phi(\al,\be)uv\hspace{3.7cm}(4.12)\end{eqnarray*}
for $u\in{\cal A}_{\al}$ and $v\in{\cal A}_{\be}$ (cf. (2.2.12), (2.2.14)).
It can be verified that the pair $({\cal A},[\cdot,\cdot])$ forms a Lie algebra. Obviously, $1$ is a central element of ${\cal A}$. Form a quotient algebra
$$H={\cal A}/\Bbb{F},\eqno(4.13)$$
whose induced Lie bracket is also denoted by $[\cdot,\cdot]$ when the context is clear.
We call the Lie algebra  $(H,[\cdot,\cdot])$ a {\it  generalized Lie algebra of Hamiltonian type}.
\psp

The following fact in linear algebra will be used.
\psp

{\bf Lemma 4.1}. 
 {\it Let} $T$ {\it be a linear transformation on a vector space} $U$ {\it and let} $U_1$ {\it be a subspace of} $U$ {\it such that} $T(U_1)\subset U_1$. {\it Suppose that} $u_1,u_2,...,u_n$ {\it are eigenvectors of} $T$ {\it corresponding to different eigenvalues. If}  $\sum_{p=1}^nu_p\in U_1$, {\it then} $u_1,u_2,...,u_n\in U_1$.
\psp

{\bf Theorem 4.2}. {\it The quotient algebra} $(H,[\cdot,\cdot])$ {\it is a simple Lie algebra if}  $\vec{\cal J}\neq \{\vec{0}\}$ {\it or} $\sgm=0$. {\it If} $\vec{\cal J}=\{\vec{0}\}$ {\it and} $\sgm\neq 0$, {\it then} $H^{(1)}=[H,H]$ {\it is a simple Lie algebra and} $H=H^{(1)}\oplus(\Bbb{F}x_1^{\sgm}+\Bbb{F})$.
\psp

{\it Proof}. We divide our proof as two parts. 
\psp

{\it Proof of the First Statement in the Theorem}
\psp
 
Note that the first statement in the theorem is equivalent to that any ideal of ${\cal A}$ that strictly contains $\Bbb{F}$ is equal to ${\cal A}$. Let ${\cal I}$ be an ideal of ${\cal A}$ such that ${\cal I}\supset \Bbb{F}$ and ${\cal I}\neq \Bbb{F}$. 
\psp

{\it Step 1}. $x_1^{\al}\in {\cal I}$ for some $0\neq \al\in\G$ or $x_2^{1_{[p]}}\in {\cal I}$ for some ${\cal J}_p=\Bbb{N}$. 
\psp

For any $\vec{j}\in\vec{\cal J}$, we define
$$|\vec{j}|=\sum_{p=1}^nj_p.\eqno(4.14)$$
Set 
$${\cal A}_k=\mbox{span}\:\{x^{\al,\vec{j}}\mid (\al,\vec{j})\in\G\times\vec{\cal J},\;|\vec{j}|\leq k\}\qquad\for\;\;k\in\Bbb{N}.\eqno(4.15)$$
For convenience, we let
$${\cal A}_{-1}=\emptyset.\eqno(4.16)$$
Moreover, we define 
$$\hat{k}=\mbox{min}\;\{k\in\Bbb{N}\mid ({\cal A}_k\bigcap{\cal I})\setminus \Bbb{F}\neq \emptyset\}.\eqno(4.17)$$
For any $u\in ({\cal A}_{\hat{k}}\bigcap{\cal I})\setminus \Bbb{F}$, we write:
$$u=u_{ld}+\td{u}\qquad\mbox{with}\;\;\td{u}\in {\cal A}_{\hat{k}-1}+\Bbb{F}\eqno(4.18)$$
and
$$u_{ld}=\sum_{(0,\vec{0})\neq (\al,\vec{j})\in\G\times\vec{\cal J},\:|\vec{j}|=\hat{k}}a_{\al,\vec{j}}x^{\al,\vec{j}},\qquad\;\;a_{\al,\vec{j}}\in\Bbb{F}\eqno(4.19)$$
and  define
$$\flat(u)=|\{\al\in\G\mid a_{\al,\vec{j}}\neq 0\;\mbox{for some}\;\vec{j}\in\vec{\cal J},\;|\vec{j}|=\hat{k}\}|.\eqno(4.20)$$
Furthermore, we set
$$\flat=\min\{\flat(v)\mid v\in ({\cal A}_{\hat{k}}\bigcap{\cal I})\setminus \Bbb{F}\}.\eqno(4.21)$$

Let $u\in ({\cal A}_{\hat{k}}\bigcap{\cal I})\setminus \Bbb{F}$ such that $\flat(u)=\flat$. Write $u$ as in (4.18) and (4.19). We set
$$q'=m+q,\;\;(m+q)'=q,\;\;\es_q=1,\;\;\es_{m+q}=-1\qquad\for\;\;q\in\ol{1,m}.\eqno(4.22)$$
If $\vf_p\neq 0$ for some $p\in\ol{m_1+1,m}\bigcup \ol{m+m_1+1,n}$, then ${\cal J}_{p'}=\Bbb{N}$ by (4.4). Thus
$$[x_2^{1_{[p']}},u]\equiv \sum_{(0,\vec{0})\neq (\al,\vec{j})\in\G\times\vec{\cal J},\:|\vec{j}|=\hat{k}}\es_{p'}a_{\al,\vec{j}}\vf_p(\al)x^{\al,\vec{j}}\;\;(\mbox{mod}\;{\cal A}_{\hat{k}-1}).\eqno(4.23)$$
By the minimality of $\flat(u)$ (cf. (4.21)) and Lemma 4.1,
$$\vf_p(\al)=\vf_p(\be)\;\;\;\mbox{whenever}\;\;a_{\al,\vec{j}}a_{\be,\vec{j}'}\neq 0.\eqno(4.24)$$

Assume that  $\vf_p\neq 0$ and $\vf_{p'}\neq 0$ for some $p\in\ol{1,m_1}\bigcup\ol{m+1,m+m_1}$.  Then for any $\tau\in\rad_{\phi}\bigcap (\bigcap_{p\neq q\in\ol{1,n}}\kn_{\vf_q})$, we have
\begin{eqnarray*}& &[x_1^{-\tau-2\sgm_p},[u,x_1^{\tau}]]\\&\equiv& \sum_{(0,\vec{0})\neq (\al,\vec{j})\in\G\times\vec{\cal J},\:|\vec{j}|=\hat{k}}\vf_p(\tau)\vf_{p'}(\al)(\vf_p(\tau)\vf_{p'}(\al-\sgm_p)+2\vf_p(\sgm_p)\vf_{p'}(\al+\sgm_p)\\& &-2\vf_p(\al+\sgm_p)\vf_{p'}(\sgm_p))a_{\al,\vec{j}}x^{\al,\vec{j}}\;\;(\mbox{mod}\:{\cal A}_{\hat{k}-1}).\hspace{6.3cm}(4.25)\end{eqnarray*}
Since $\vf_p(\tau)$ takes an infinite number of elements in $\Bbb{F}$ if $\tau$ varies in $\rad_{\phi}\bigcap (\bigcap_{p\neq q\in\ol{1,n}}\kn_{\vf_q})$ by (4.7),
the coefficients of $\vf_p(\tau)^2$ in (4.25) show
$$\vf_{p'}(\al)\vf_{p'}(\al-\sgm_p)=\vf_{p'}(\be)\vf_{p'}(\be-\sgm_p)\;\;\;\mbox{whenever}\;\;a_{\al,\vec{j}}a_{\be,\vec{j}'}\neq 0\eqno(4.26)$$
by the minimality of $\flat(u)$ (cf. (4.21)) and Lemma 4.1.
Moreover, (4.26) is equivalent to
$$\vf_{p'}(\al)=\vf_{p'}(\be)\;\;\mbox{or}\;\;\vf_{p'}(\al+\be-\sgm_p)=0\;\;\;\mbox{whenever}\;\;a_{\al,\vec{j}}a_{\be,\vec{j}'}\neq 0.\eqno(4.27)$$
Assume that there exist $\al,\be\in\G$ such that $\vf_{p'}(\al)\neq \vf_{p'}(\be),\;\vf_{p'}(\al+\be-\sgm_p)=0$ and $a_{\al,\vec{j}}a_{\be,\vec{j}'}\neq 0$. We may assume that $\vf_{p'}(\al)\neq 0$. Since $\vf_p\neq 0$ and $\vf_{p'}\neq 0$, we can choose 
$$\tau\in\rad_{\phi}\bigcap (\bigcap_{p\neq q\in\ol{1,n}}\kn_{\vf_q})\setminus\kn_{\vf_p},\;\;\;\tau'\in\rad_{\phi}\bigcap (\bigcap_{p'\neq q\in\ol{1,n}}\kn_{\vf_q})\setminus\kn_{p'},\eqno(4.28)$$
Such that
$$\vf_{p'}(\tau'+\al)\neq 0,\;\;\;\vf_p(\tau-\sgm_p)\vf_{p'}(\al)\neq \vf_p(\al)\vf_{p'}(\tau'-\sgm_p)\eqno(4.29)$$
by (4.7). We have
\begin{eqnarray*}& &[x_1^{\tau+\tau'-\sgm_p},u]\equiv \sum_{(0,\vec{0})\neq (\gm,\vec{l})\in\G\times\vec{\cal J},\:|\vec{l}|=\hat{k}}\es_pa_{\gm,\vec{l}}(\vf_p(\tau-\sgm_p)\vf_{p'}(\gm)\\& &-\vf_p(\gm)\vf_{p'}(\tau'-\sgm_p))x^{\gm+\tau+\tau',\vec{l}}\;\;(\mbox{mod}\;{\cal A}_{\hat{k}-1}).\hspace{6.5cm}(4.30)\end{eqnarray*}
Since $\flat(u)$ is minimum, $\al+\tau+\tau'\neq 0$ due to $\vf_{p'}(\al+\tau+\tau')\neq 0$ and
$$\es_pa_{\al,\vec{j}}(\vf_p(\tau-\sgm_p)\vf_{p'}(\al)-\vf_p(\al)\vf_{p'}(\tau'-\sgm_p))\neq 0,\eqno(4.31)$$
 we have
$$\es_pa_{\be,\vec{j}'}(\vf_p(\tau-\sgm_p)\vf_{p'}(\be)-\vf_p(\be)\vf_{p'}(\tau'-\sgm_p))\neq 0\eqno(4.32)$$
by Lemma 4.1. But 
$$\vf_{p'}(\al+\tau+\tau')\neq \vf_{p'}(\be+\tau+\tau')\eqno(4.33)$$
and
$$\vf_{p'}((\al+\tau+\tau')+(\be+\tau+\tau')-\sgm_p)=2\vf_{p'}(\tau')\neq 0,\eqno(4.34)$$
which contradicts (4.27) with $u$ replaced by $[x_1^{\tau+\tau'-\sgm_p},u]$. Thus the first equation in (4.27) holds. 

Assume $\vf_p=0$ with $p\in\ol{1,m_1}\bigcup \ol{m+1,m+m_1}$. Then ${\cal J}_p=\Bbb{N}$ by (4.2). 
If $j_p\neq 0$ for some $a_{\al,\vec{j}}\neq 0$, then $[x_1^{l\sgm_p},u]\in {\cal I}\bigcap {\cal A}_{\hat{k}-1}\setminus\Bbb{F}$ for some $0<l\in \Bbb{N}$ such that $(l+1)\sgm_p+\al\neq 0$ by (4.9), which contradicts (4.17).
Assume that $j_p=0$ whenever $a_{\al,\vec{j}}\neq 0$. In this case,
 we get 
$$[x_2^{-\sgm_p,1_{[p]}},u]\equiv \sum_{(0,\vec{0})\neq (\al,\vec{j})\in\G\times\vec{\cal J},\:|\vec{j}|=\hat{k}}a_{\al,\vec{j}}\es_p\vf_{p'}(\al)x^{\al,\vec{j}}\;\;(\mbox{mod}\:{\cal A}_{\hat{k}-1}),\eqno(4.35)$$
which implies the first equation in (4.27). 
 Therefore, we have proved that (4.24) for any $p\in\ol{1,n}$.

If $a_{\be,\vec{j}'}\neq 0$, then
$${\cal I}\ni [x_1^{\be},u]\equiv\sum_{(0,\vec{0})\neq (\al,\vec{j})\in\G\times\vec{\cal J},\:|\vec{j}|=\hat{k}}a_{\al,\vec{j}}\phi(\be,\al)x^{\al+\be,\vec{j}}\;\;(\mbox{mod}\;{\cal A}_{\hat{k}-1})\eqno(4.36)$$
by (4.12) and (4.24). Since $\flat(u)$ is minimum and $\phi(\be,\be)=0$, we have
$$\phi(\al,\be)=0\;\;\;\mbox{whenever}\;\;a_{\al,\vec{j}}a_{\be,\vec{j}'}\neq 0.\eqno(4.37)$$
Now for any $\gm\in\bigcap_{q\in\mho}\kn_{\vf_q}$, we have
$${\cal I}\ni [[x_1^{-\gm},[x_1^{\gm},u]]\equiv\sum_{(0,\vec{0})\neq (\al,\vec{j})\in\G\times\vec{\cal J},\:|\vec{j}|=\hat{k}}-a_{\al,\vec{j}}\phi(\gm,\al)^2x^{\al,\vec{j}}\;\;(\mbox{mod}\;{\cal A}_{\hat{k}-1}).\eqno(4.38)$$
Since $\flat(u)$ is minimum, by Lemma 4.1,  we get
$$\phi(\gm,\al)^2=\phi(\gm,\be)^2\;\;\;\mbox{whenever}\;\;a_{\al,\vec{j}}a_{\be,\vec{j}'}\neq 0,\eqno(4.39)$$
which is equivalent to
$$\phi(\gm,\al)=\pm \phi(\gm,\be)\;\;\;\mbox{whenever}\;\;a_{\al,\vec{j}}a_{\be,\vec{j}'}\neq 0.\eqno(4.40)$$

Assume that there exist $\al,\be\in\G$ such that $\phi(\gm,\al)=-\phi(\gm,\be)\neq 0$ and $a_{\al,\vec{j}}a_{\be,\vec{j}'}\neq 0$. Then
$${\cal I}\ni [x_1^{\gm},u]\equiv\sum_{(0,\vec{0})\neq (\al,\vec{j})\in\G\times\vec{\cal J},\:|\vec{j}|=\hat{k}}a_{\al,\vec{j}}\phi(\gm,\al)x^{\al+\gm,\vec{j}}\;\;(\mbox{mod}\;{\cal A}_{\hat{k}-1}).\eqno(4.41)$$
Thus $\flat([x_1^{\gm},u])=\flat(u)$ and 
$$a_{\al,\vec{j}}\phi(\gm,\al)a_{\be,\vec{j}}\phi(\gm,\be)\neq 0.\eqno(4.42)$$
But 
$$\phi(\be+\gm,\al+\gm)=\phi(\gm,\al)+\phi(\be,\gm)=2\phi(\gm,\al)\neq 0\eqno(4.43)$$
by (4.37) and the skew-symmetry of $\phi$. Equation (4.43) contradicts (4.37) if we replace $u$ by $[x_1^{\gm},u]$. Hence we have
$$\phi(\gm,\al)=\phi(\gm,\be)\;\;\;\mbox{whenever}\;\;a_{\al,\vec{j}}a_{\be,\vec{j}'}\neq 0\eqno(4.44)$$
for any $\gm\in\bigcap_{q\in\mho}\kn_{\vf_q}$. By (4.8), (4.24) and (4.44), we obtain
$$\al=\be\;\;\;\mbox{whenever}\;\;a_{\al,\vec{j}}a_{\be,\vec{j}'}\neq 0.\eqno(4.45)$$

Let $a_{\al,\vec{j}}\neq 0$ be fixed. If $\hat{k}=0$, then we have $x_1^{\al}\in {\cal I}$. Assume that $\hat{k}>0$. 
\pse

{\it Case 1}. $\al=0$ and $\hat{k}=1$.
\pse

In this case, 
$$u=\sum_{p=1}^{n}a_px_2^{1_{[p]}}+\sum_{\be\in\G}b_{\be}x_1^{\be}.\eqno(4.46)$$
For convenience of stating things, we set
$$\sgm_p=\sgm_{m+p}=0\qquad\for\;\;p\in\ol{m_1+1,m}.\eqno(4.47)$$
If $\vf_p\neq 0$ with $p\in\ol{1,n}$, we choose  $-\sgm_p\neq \tau\in (\bigcap_{p\neq q\in\ol{1,n}}\kn_{\vf_q})\setminus\kn_{\vf_p}$ by (4.6), (4.7) and have
$$[x_1^{\tau},u]=\es_p\vf_p(\tau)a_{p'}x_1^{\tau+\sgm_p}+\sum_{\be\in\G}b_{\be}(\phi(\tau,\be)x_1^{\be+\tau}+\es_p\vf_p(\tau)\vf_{p'}(\be)x_1^{\be+\tau+\sgm_p})\in {\cal I}.\eqno(4.48)$$
If $a_{p'}\neq 0$, then $[x_1^{\tau},u]\not\in\Bbb{F}$, which contradicts the assumption $\hat{k}=1$. 

If $\vf_p=0$, then ${\cal J}_p=\Bbb{N}$. We have
$$[x_2^{1_{[p]}},u]=\es_pa_{p'}x_1^{\sgm_p}+\sum_{\be\in\G}\es_pb_{\be}\vf_{p'}(\be)x_1^{\sgm_p+\be}\in {\cal I}.\eqno(4.49)$$
If $a_{p'}\neq 0$ and $\es_pb_{\be}\vf_{p'}(\be)\neq 0$ for some $\be\in\G$, we get a contradiction to the assumption $\hat{k}=1$. We assume $a_{p'}\neq 0$ for some $p'\in\ol{1,n}$, then we have $\vf_p=0$ by (4.48) and
$$\es_pb_{\be}\vf_{p'}(\be)=0\qquad\for\;\;\be\in\G\eqno(4.50)$$
by (4.49). Note ${\cal J}_p=\Bbb{N}$ in this case. So
$$[x_2^{2_{[p]}},u]=2\es_pa_{p'}x_2^{1_{[p]}}\in{\cal I}.\eqno(4.51)$$
Thus the conclusion of Step 1 holds.
\pse

{\it Case 2}. $\al\neq 0$ or $\hat{k}>1$.  $j_p'\neq 0$ for some  $a_{\al,\vec{j}'}\neq 0$ with $p\in\ol{1,n}$ and $\vf_p(\al)=0$.
\pse

If $p\in\ol{1,m_1}\bigcap\ol{m+1,m+m_1}$ and $\vf_{p'}\neq 0$, we choose $\tau'\in \rad_{\phi}\bigcap (\bigcap_{p'\neq q\in\ol{1,n}}\kn_{\vf_q})\setminus\kn_{\vf_{p'}}$ such that $\vf_{p'}(\al+\tau'+\sgm_p)\neq 0$ by (4.7). Then 
$$ [u,x_1^{\tau'}]\in ({\cal I}\bigcap{\cal A}_{\hat{k}-1})\setminus \Bbb{F},\eqno(4.52)$$
which contradicts  (4.17).

 If $p\in\ol{m_1+1,m}\bigcap\ol{m+m_1+1,n},\;\vf_{p'}\neq 0$, we let $\tau'\in (\bigcap_{p'\neq q\in\ol{1,n}}\kn_{\vf_q})\setminus\kn_{\vf_{p'}}$ such that $\vf_{p'}(\al+\tau')\neq 0$ by (4.6). If $\phi(\tau',\al)=0$, then we have (4.52). Otherwise,
\begin{eqnarray*}v&=&[x_1^{\al+\tau'},[x_1^{\tau'},u]]\\&\equiv& \sum_{\vec{j}''\in\vec{\cal J},\;|\vec{j}''|=\hat{k}}\phi(\tau',\al)a_{\al,\vec{j}''}(\vf_{p'}(\tau')+\vf_{p'}(\al))j_p''x^{\al,\vec{j}''-1_{[p]}}+w\;\;(\mbox{mod}\;{\cal A}_{\hat{k}-2}),\hspace{1.3cm}(4.53)\end{eqnarray*}
where
$$w=\sum_{(\be,\vec{j}'')\in\G\times \vec{\cal J},\;|\vec{j}''|=\hat{k}-1}b_{\be,\vec{j}''}x^{\be,\vec{j}''}\eqno(4.54)$$
such that $b_{\be,\vec{j}''}\in\Bbb{F}$ are independent of $\vf_{p'}(\tau')$ and the number of nonzero $b_{\be,\vec{j}''}$ is bounded with respect to $\vf_{p'}(\tau')$. Since $|\vf_{p'}(\Bbb{Z}\tau')|=\infty$, $\vf_{p'}(\tau')$ takes an infinite number of elements in $\Bbb{F}$ when $\tau'$ varies in $\tau'\in (\bigcap_{p'\neq q\in\ol{1,n}}\kn_{\vf_q})\setminus\kn_{\vf_{p'}}$ such that $\phi(\tau',\al)\neq 0$. Hence we can choose $\tau'$ such  that $v\in {\cal A}_{\hat{k}-1}\setminus \Bbb{F}$, which contradicts (4.17).

Assume $\vf_{p'}=0$. We have ${\cal J}_{p'}=\Bbb{N}$ by (4.2).  Then
$$ [u,x_2^{1_{[p']}}]\in ({\cal I}\bigcap{\cal A}_{\hat{k}-1})\setminus \Bbb{F},\eqno(4.55)$$
which contradicts (4.17). 
\pse

{\it Case 3}. $j_p'\neq 0$ for some  $a_{\al,\vec{j}'}\neq 0$ with $p\in\ol{m_1+1,m}\bigcup\ol{m+m_1+1,n}$ and $\vf_p(\al)\neq 0$.
\pse

Note that $[x_1^{-\al},u]\in{\cal I}\bigcap {\cal A}_{\hat{k}-1}$. By (4.17), we have
$$[x_1^{-\al},u]=\lmd\in\Bbb{F}.\eqno(4.56)$$
 Thus we get
\begin{eqnarray*}& &{\cal I}\ni[u,x^{-\al}_1x_2^{1_{[p']}}]=[u,x^{-\al}_1]x_2^{1_{[p']}}+x_1^{-\al}[u,x_2^{1_{[p']}}]\\&\equiv&\lmd x_2^{1_{[p']}}+\sum_{\vec{j}''\in\vec{\cal J},\;|\vec{j}''|=\hat{k}}\vf_p(\al)\es_pa_{\al,\vec{j}''}x^{0,\vec{j}''}\;\;(\mbox{mod}\;{\cal A}_{\hat{k}-1}).\hspace{4.9cm}(4.57)\end{eqnarray*}
Replacing $u$ by $[u,x^{-\al}_1x_2^{1_{[p']}}]$, we go back to Case 2 if $\hat{k}>1$ and to Case 1 if $\hat{k}=1$ because $j'_p\neq 0$.
\pse

{\it Case 4}.  There exists $p\in\ol{1,m_1}\bigcup\ol{m+1,m+m_1}$ such that $\vf_p(\al)\neq 0,\;\vf_{p'}\neq 0$ and $j_p>0$.

We can choose $\tau' \in \rad_{\phi}\bigcap (\bigcap_{p'\neq q\in\ol{1,n}}\kn_{\vf_q})\setminus\kn_{\vf_{p'}}$ such that $\vf_{p'}(\tau'+\sgm_p)\vf_p(\al)\neq \vf_p(\sgm_p)\vf_{p'}(\al)$ by (4.7). Note that
\begin{eqnarray*}{\cal I}\ni [u,x_1^{-\al-\tau'-\sgm_p}]&\equiv& \es_{p'}(\vf_{p'}(\tau'+\sgm_p)\vf_p(\al)-\vf_p(\sgm_p)\vf_{p'}(\al))\\& &\sum_{\vec{j}'\in\vec{\cal J},\:|\vec{j}'|=\hat{k}}a_{\al,\vec{j}'}x^{\tau',\vec{j}'}\;\;(\mbox{mod}\;{\cal A}_{\hat{k}-1}).\hspace{4.8cm}(4.58)\end{eqnarray*}
Since $\vf_p(\tau')=0$ and $j_p>0$, we go back to Case 2 with $u$ replaced by $[u,x_1^{-\al-\tau'-\sgm_p}]$.
\pse

{\it Case 5}. There exists $p\in\ol{1,m_1}\bigcup\ol{m+1,m+m_1}$ such that $\vf_p(\al)\neq 0,\;\vf_{p'}=0,\;j_p>0$ and $j_{p'}=0$ for all $a_{\al,\vec{j}'}\neq 0$.
\pse

By (4.2), we have ${\cal J}_{p'}=\Bbb{N}$. Note that (4.56) holds. As (4.57), we have
\begin{eqnarray*}& &{\cal I}\ni[u,x^{-\al}_1x^{-\sgm_p,1_{[p']}}]=[u,x^{-\al}_1]x^{-\sgm_p,1_{[p']}}+x_1^{-\al}[u,x^{-\sgm_p,1_{[p']}}]\\&\equiv&\lmd x^{-\sgm_p,1_{[p']}}+\sum_{\vec{j}''\in\vec{\cal J},\;|\vec{j}''|=\hat{k}}\vf_p(\al)\es_pa_{\al,\vec{j}''}x^{0,\vec{j}''}\;\;(\mbox{mod}\;{\cal A}_{\hat{k}-1}).\hspace{4.3cm}(4.59)\end{eqnarray*}
If $\hat{k}>1$ or $\lmd=0$, we go back to Case 2 if we replace $u$ by $[u,x^{-\al}_1x^{-\sgm_p,1_{[p']}}]$. Assume that $\lmd\neq 0$ and $\hat{k}$=1. Then we have
$${\cal I}\ni [[u,x^{-\al}_1x^{-\sgm_p,1_{[p']}}],x_2^{1_{[p']}}]\equiv-\es_p\lmd\vf_p(\sgm_p)x_2^{1_{[p']}}\;\;(\mbox{mod}\;{\cal A}_0).\eqno(4.60)$$
Replacing $u$ by $[[u,x^{-\al}_1x^{-\sgm_p,1_{[p']}}],x_2^{1_{[p']}}]$, we go back to Case 1.
 
\psp

This completes the proof of the conclusion in Step 1. 
\psp

{\it Step 2}. The conclusion of Step 1 implies ${\cal I}={\cal A}$.
\psp

{\it Case 1}. $x_2^{1_{[p]}}\in {\cal I}$ for some $p\in\ol{1,n}$.
\pse

In this case, ${\cal J}_p=\Bbb{N}$. For any $(\be,\vec{j}')\in\G\times\vec{\cal J}$, we have
$$[x_2^{1_{[p]}},x^{\be,\vec{j}'}]=\es_p (\vf_{p'}(\be)x^{\be+\sgm_p,\vec{j}'}+j_{p'}'x^{\be+\sgm_p,\vec{j}'-1_{[p']}})\in {\cal I}.\eqno(4.61)$$
If ${\cal J}_{p'}=\Bbb{N}$, then by (4.61) and induction on $j_{p'}'$, we can prove ${\cal I}={\cal A}.$
Assume ${\cal J}_{p'}=\{0\}$. By (4.2), $\vf_{p'}\neq 0$. We choose $\tau'\in (\bigcap_{p'\neq q\in\ol{1,n}}\kn_{\vf_q})\setminus\kn_{\vf_{p'}}$ such that $\vf_{p'}(\tau'-\sgm_p)\neq 0$ by (4.6) and (4.7). Moreover, we can assume $\tau'\in\rad_{\phi}$ if $p'\in\ol{1,m_1}\bigcup\ol{m+1,m+m_1}$ by (4.7). We have
$$[x_2^{1_{[p]}},x_1^{\tau'-\sgm_p}]=\es_p\vf_{p'}(\tau'-\sgm_p)x_1^{\tau'}\in {\cal I}.\eqno(4.62)$$
So $x_1^{\tau'}\in {\cal I}$. Moreover, for any $(\be,\vec{j}')\in\G\times\vec{\cal J}$,
$$[x_1^{\tau'},x^{\be,\vec{j}'}]=\es_{p'} \vf_{p'}(\tau')(\vf_p(\be)x^{\be+\tau'+\sgm_p,\vec{j}'}+j_p'x^{\be+\sgm_p+\tau',\vec{j}'-1_{[p]}})\in {\cal I}\eqno(4.63)$$
if $p\in\ol{1,m_1}\bigcup \ol{m+1,m+m_1}$ and
$$[x_1^{\tau'},x^{\be,\vec{j}'}]=\phi(\tau',\be)x^{\be+\tau',\vec{j}'}+\es_{p'} \vf_{p'}(\tau')j'_px^{\be+\tau',\vec{j}'-1_{[p]}}\in {\cal I}\eqno(4.64)$$
if $p\in\ol{m_1+1,m}\bigcup \ol{m+m_1+1,n}$. Thus by (4.63), (4.64) and induction on $j'_p$, we can prove ${\cal I}={\cal A}.$
\pse

{\it Case 2}. $x_1^{\al}\in {\cal I}$ for some $0\neq\al\in \G$ and $\vf_p(\al)\neq 0$ for some $p\in\ol{1,n}$.
\pse

Assume $p\in\ol{1,m_1}\bigcup \ol{m+1,m+m_1}$ and $\vf_{p'}=0$ or $p\in\ol{m_1+1,m}\bigcup \ol{m+m_1+1,n}$. We have ${\cal J}_{p'}=\Bbb{N}$ by (4.2) and (4.4). Moreover,
$$[x_1^{\al},x^{-\al-\sgm_p,2_{[p']}}]=2\es_p\vf_p(\al)x_2^{1_{[p']}}\in {\cal I}.\eqno(4.65)$$
So $x_2^{1_{[p']}}\in {\cal I}$. Hence ${\cal I}={\cal A}$ by Case 1. Next we consider $p\in\ol{1,m_1}\bigcup \ol{m+1,m+m_1}$ and $\vf_{p'}\neq 0$. By (4.58) with $u=x_1^{\al}$, we obtain $x_1^{\tau'}\in {\cal I}$ for some
$\tau'\in\rad_{\phi}\bigcap (\bigcap_{p'\neq q\in\ol{1,n}}\kn_{\vf_q})\setminus\kn_{\vf_{p'}}$. Furthermore,  by (4.63), we have ${\cal I}={\cal A}$ if ${\cal J}_p=\Bbb{N}$ and otherwise,
$$x^{\be,\vec{j}'}\in {\cal I}\qquad\mbox{for any}\;\;(\be,\vec{j}')\in\G\times\vec{\cal J},\;\vf_p(\be)\neq\vf_p(\sgm_p).\eqno(4.66)$$
Choose any $\tau\in (\bigcap_{p\neq q\in\ol{1,n}}\kn_{\vf_q})\setminus\kn_{\vf_p}$ such that $\vf_p(\tau)\neq \vf_p(\sgm_p)$. We have $x_1^{\tau+\tau'}\in {\cal I}$ by (4.66) and $\vf_{p'}(\tau+\tau')=\vf_{p'}(\tau')\neq 0$. Exchanging positions of $p$ and $p'$, we have the similar conclusion. If $\vf_q\neq 0$ for some $q\in\ol{1,n}\setminus\{p,p'\}$, then we choose any $\tau_q\in(\bigcap_{q\neq r\in\ol{1,n}}\kn_{\vf_r})\setminus\kn_{\vf_q}$ and have $x_1^{\tau+\tau_q}\in {\cal I},\;\vf_q(\tau+\tau_q)=\vf_q(\tau_q)\neq 0$. Thus we have the same conclusion for $(q,q')$ as that for $(p,p')$. Assume that $\vf_q=\vf_{q'}=0$ for some $q\in\ol{1,n}\setminus\{p,p'\}$, then ${\cal J}_q={\cal J}_{q'}=\Bbb{N}$. Since $x^{\tau,1_{[q]}}\in {\cal I}$, we have
$$[x^{\tau,1_{[q]}},x^{-\tau,2_{[q']}}]=2\es_qx_2^{1_{[q']}}\in {\cal I},\eqno(4.67)$$
by which and Case 1, we get ${\cal I}={\cal A}$. In summary, we have ${\cal I}={\cal A}$ if $\vec{\cal J}\neq\{\vec{0}\}$ and
$$x_1^{\be}\in {\cal I}\qquad\for\;\;\be\in\G\;\mbox{such that}\;\be-\sgm\not\in \bigcap_{r=1}^{n}\kn_{\vf_r}\eqno(4.68)$$
(cf. (4.10)) if $\vec{\cal J}=\{\vec{0}\}$ , where $m_1=m$ by (4.2) and (4.4).

Assume $\vec{\cal J}=\{\vec{0}\}$. Let $\be\in\G$ such that $0\neq \be-\sgm\in \bigcap_{r=1}^{n}\kn_{\vf_r}$. By (4.8), there exists $\gm\in \bigcap_{r=1}^{n}\kn_{\vf_r}$ such that $\phi(\gm,\be)=\phi(\gm,\be-\sgm)\neq 0$, where we have used the fact that $\sgm\in \rad_{\phi}$ by (4.9) and (4.10). Moreover, we choose $\tau \in\rad_{\phi}\bigcap (\bigcap_{p\neq q\in\ol{1,n}}\kn_{\vf_q})\setminus\kn_{\vf_p}$. Since $\vf_p(\be-\tau-\gm)=\vf_p(\sgm_p)-\vf(\tau)\neq \vf_p(\sgm_p)$, we have $x_1^{\be-\tau-\gm}\in {\cal I}$ by (4.66). Note that
$$[x_1^{\gm+\tau},x_1^{\be-\gm-\tau}]=\phi(\gm,\be)x_1^{\be}+\es_p\vf_p(\tau)\vf_{p'}(\sgm_p)x^{\be+\sgm_p}\in {\cal I},\eqno(4.69)$$
where we have used the fact that $\vf_{p'}(\be)=\vf_{p'}(\sgm_p)$ and $\vf_{p'}(\tau)=0$. Moreover, we have $\es_p\vf_p(\tau)\vf_{p'}(\sgm_p)x^{\be+\sgm_p}\in {\cal I}$ by (4.68). Hence $x_1^{\be}\in {\cal I}$. Therefore, we obtain
$$x_1^{\be}\in {\cal I}\qquad\for\;\;\sgm\neq\be\in\G. \eqno(4.70)$$
\pse

{\it Case 3}. $x_1^{\al}\in {\cal I}$ for some $0\neq\al\in \G$ and $\vf_p(\al)=0$ for any $p\in\ol{1,n}$.
\pse

By (4.8), there exists $\be\in\G$ such that $\phi(\be,\al)\neq 0$. If $\vf_p\neq 0$ for some $p\in\ol{1,n}$, we can choose $\tau\in(\bigcap_{p\neq q\in\ol{1,n}}\kn_{\vf_q})\setminus\kn_{\vf_p}$ such that $\vf_p(\tau+\al+\be)\neq 0$ and $\phi(\be+\tau,\al)\neq 0$. Since
$$\phi(\be+\tau,\al)^{-1}[x_1^{\tau+\be},x_1^{\al}]=x_1^{\tau+\al+\be}\in {\cal I},\eqno(4.71)$$
we go back to Case 2. Assume that $\vf_p=0$ for any $p\in\ol{1,n}$ and $n>0$. By (4.2), we have  ${\cal J}_1=\Bbb{N}$. Moreover,
$$[[x^{\be,1_{[1]}},x_1^{\al}],x^{-\al-\be,2_{[m+1]}}]=2\phi(\be,\al)x_2^{1_{[m+1]}}\in {\cal I},\eqno(4.72)$$
by which and Case 1, ${\cal I}={\cal A}$. Now we assume $n=0$. In this situation, (4.8) becomes $\rad_{\phi}=\{0\}$.
Note
$$[x_1^{\be-\al},x_1^{\al}]=\phi(\be,\al)x_1^{\be}\in {\cal I}\qquad\mbox{for any}\;\;\be\in\G.\eqno(4.73)$$
This shows that
$$\begin{array}{l}x_1^{\gm}\in {\cal I}\;\for\;\gm\in\G,\;\phi(\gm,\al)\neq 0\;\mbox{or there exists}\\\be\in\G\;\mbox{such that}\;\phi(\gm,\be)\phi(\be,\al)\neq 0.\end{array}\eqno(4.74)$$
If (4.74) does not imply ${\cal I}={\cal A}$, then there exists $0\neq \gm\in \G$ such that $\phi(\gm,\al)=0$ and
$$\phi(\gm,\be)\phi(\be,\al)=0\qquad\mbox{for any}\;\;\be\in \G.\eqno(4.75)$$
Since $\rad_{\phi}=\{0\}$, there exists $\al_0,\gm_0\in \G$ such that $\phi(\al_0,\al)\neq 0$ and $\phi(\gm,\gm_0)\neq 0$. By (4.75), $\phi(\gm_0,\al)=\phi(\gm,\al_0)=0$. However,
$$\phi(\gm,\al_0+\gm_0)\phi(\al_0+\gm_0,\al)=\phi(\gm,\gm_0)\phi(\al_0,\al)\neq 0,\eqno(4.76)$$
which contradicts (4.75) with $\be=\al_0+\gm_0$. Thus ${\cal I}={\cal A}$ if $n=0$. This completes the proof of the first statement in the theorem.
\psp

{\it Proof of the Second Statement in the Theorem}
\psp

Now $\vec{\cal J}=\{\vec{0}\}$ and $\sgm \neq 0$. By (4.4), we have $m_1=m$. 
Set
$$\hat{H}=\mbox{span}\{x_1^{\al}\mid \sgm\neq \al\in\G\}.\eqno(4.77)$$
For any $\al,\be\in\G$, we have:
$$[x_1^{\al},x_1^{\be}]=\phi(\al,\be)x_1^{\al+\be}+\sum_{p=1}^{m}(\vf_p(\al)\vf_{m+p}(\be)-\vf_{m+p}(\al)\vf_p(\be))x_1^{\al+\be+\sgm_p}.\eqno(4.78)$$
If $\al+\be=\sgm$, then $\phi(\al,\be)=\phi(\al,\sgm-\al)=0$ by (4.9) and (4.10). For $p\in\ol{1,m}$, if $\al+\be+\sgm_p=\sgm$, then $\be=\sum_{p\neq q\in\ol{1,m}}\sgm_q-\al$. So $\vf_p(\be)=-\vf_p(\al)$ and $\vf_{m+p}(\be)=-\vf_{m+p}(\al)$, which implies $\vf_p(\al)\vf_{m+p}(\be)-\vf_{m+p}(\al)\vf_p(\be)=0$.
Thus we have $[{\cal A},{\cal A}]\subset \hat{H}$. On the other hand, $[{\cal A},{\cal A}]\supset \hat{H}$ by (4.7), (4.8), (4.63), (4.69) and (4.73)-(4.76). Hence,
$$\hat{H}=[{\cal A},{\cal A}].\eqno(4.79)$$
Replacing ${\cal A}$ by $\hat{H}$ in the proof of the first statement, we obtain (4.70), which implies that 
$H^{(1)}=\hat{H}/\Bbb{F}$ is simple.$\qquad\Box$
\psp

{\bf Remark 4.3}.  Some special cases of the Lie algebra $H$ with $m_1=m=1$ and $\phi=0$ were studied in [X3]. 

\psp

{\bf Example}. Let $k\in\Bbb{N}$. Take any nondegenerate skew-symmetric bilinear form $\phi'$ on $\Bbb{F}^k$, where we treat $\Bbb{F}^0=\{0\}$ with $\phi'=0$ for convenience. Suppose that we have picked (2.5) so that
$${\cal J}_q=\Bbb{N}\;\;\mbox{or}\;\;{\cal J}_{q'}=\Bbb{N}\qquad\for\;\;q\in\ol{m_1+1,m}\bigcup \ol{m+m_1+1,n}\eqno(4.80)$$
(cf. (4.22)).  Let $k_1$ be the number of ${\cal J}_p =\{0\}$ with $p\in\ol{1,m_1}\bigcup \ol{m+1,m+m_1}$ and let $k_2$ be the number of ${\cal J}_q=\Bbb{N}$ with $q\in\ol{m_1+1,m}\bigcup \ol{m+m_1+1,n}$. Set $s=k_1+n-2m_1-k_2$. Pich an integer $\ell$ such that $s\leq \ell\leq 2m_1+k_2$. We define
$$\zeta_p(\al_1,...,\al_{\ell})=\al_p\qquad\for\;\;p\in\ol{1,\ell},\;(\al_1,\al_2,...,\al_{\ell})\in\Bbb{F}^{\ell}.\eqno(4.81)$$
Moreover, we extend $\phi'$ and $\zeta_p$ to $\Bbb{F}^{k+\ell}$ by
$$\phi'((\vec{\al}_1,\vec{\al}_2),(\vec{\be}_1,\vec{\be}_2))=\phi'(\vec{\al}_1,\vec{\be}_1),\;\;\zeta_p(\vec{\al}_1,\vec{\al}_2)=\zeta_p(\vec{\al}_2)\eqno(4.82)$$
for $\vec{\al_1},\vec{\be}_1\in\Bbb{F}^k,\;\vec{\al}_2,\vec{\be}_2\in\Bbb{F}^{\ell}$ and $p\in\ol{1,\ell}$. Furthermore, we define 
$$\zeta_q\equiv 0 \qquad\for\;\;q\in\ol{\ell+1,n}.\eqno(4.83)$$
Take $\G$ to be an additive subgroup of $\Bbb{F}^{k+\ell}$ containing $\Bbb{Z}^{k+\ell}$. Take any permutation $\iota$ on $\ol{1,n}$ such that
$$\iota(p)\leq\ell\;\;\mbox{if}\;\;{\cal J}_p=\{0\}\;\;\mbox{with}\;\;p\in\ol{1,m_1}\bigcup \ol{m+1,m+m_1},\eqno(4.84)$$
$$\iota(q)>\ell,\;\iota_{q'}\leq \ell\;\;\mbox{if}\;\;{\cal J}_{q'}=\{0\}\;\;\mbox{with}\;\;q\in\ol{m_1+1,m}\bigcup \ol{m+m_1+1,n}\eqno(4.85)$$
(cf. (4.22)). We let
$$\phi=\phi'(|_{\G},|_{\G}),\;\;\vf_p=\zeta_{\iota(p)}|_{\G}\qquad\for\;\;p\in\ol{1,n}.\eqno(4.86)$$
Then (4.2)-(4.4) and (4.6)-(4.8) hold. 

In particular, we can take $k=0$, $m_1=m$, (2.34) and (2.35).

\section{Algebras of Type K}

In this section, we shall construct and prove a new class of  generalized simple Lie algebras of Contact type. We shall use the notions in Section 2.  All the  notations and assumptions except (2.6) are the same as in Section 2.

Now we assume 
$$n=2m+1\;\;\mbox{for some}\;\;0<m\in\Bbb{N};\qquad \vf_p\not\equiv 0\;\;\mbox{or}\;\;{\cal J}_p=\Bbb{N}\;\;\for\;\;p\in\ol{1,2m+1}\eqno(5.1)$$
and
$$\G=\G_1+\G_2,\qquad \G_1=\kn_{\vf_n},\;\;\G_2=\bigcap_{p=1}^{2m}\kn_{\vf_p}.\eqno(5.2)$$
By (2.7), the sum $\G_1+\G_2$ is a direct sum, $\vf_n:\G_2\rightarrow \vf_n(\G_2)=\vf_n(\G)$ is
an additive group homomorphism and
$$\bigcap_{p=1}^{2m}\kn_{\vf_p|_{\G_1}}=\{0\}.\eqno(5.3)$$
Let
$$\mho_1=\{p\in\ol{1,2m}\mid \vf_p\not\equiv 0\},\;\;\mho_2=\ol{1,2m}\setminus\mho_1.\eqno(5.4)$$
Condition (2.6) is replaced by
$$\bigcap_{q\neq p\in\ol{1,2m}}\kn_{\vf_p|_{\G_1}}\setminus\kn_{\vf_q|_{\G_1}}\neq\emptyset\qquad\for\;\;q\in\mho_1.\eqno(5.5)$$
Now we choose fixed elements
$$\sgm_q=\sgm_{q'}\in \bigcap_{q,q'\neq p\in\ol{1,2m}}\kn_{\vf_p|_{\G_1}}\setminus(\kn_{\vf_q|_{\G_1}}\bigcup\kn_{\vf_{q'}|_{\G_1}})\qquad\for\;\;q\in\mho_1\eqno(5.6)$$
(cf. (4.22)) and 
$$\sgm_n\in\G_2.\eqno(5.7)$$
Up to equivalence of the following construction of the Lie algebras of type K, we can assume
$$\vf_p(\sgm_p)=-1\qquad\for\;\;p\in\mho_1.\eqno(5.8)$$
For convenience, we let
$$\sgm_p=0\qquad\for\;\;p\in\mho_2.\eqno(5.9)$$

Note by the assumption (5.1),
$${\cal J}_q=\Bbb{N}\qquad\for\;\;q\in\mho_2.\eqno(5.10)$$
 We define an operator $\ptl$ on ${\cal A}$ by
$$\ptl(x^{\al,\vec{i}})=(\sum_{p\in\mho_1}\vf_p(\al)+\sum_{q\in\mho_2}i_q)x^{\al,\vec{i}}\qquad\for\;\;(\al,\vec{i})\in\G\times\vec{\cal J}.\eqno(5.11)$$
Then $\ptl$ is a derivation of $({\cal A},\cdot)$.
Moreover, we define an algebraic operation $[\cdot,\cdot]_K$ on ${\cal A}$ by 
$$[u,v]_K=\sum_{p=1}^mx_1^{\sgm_p}(\ptl_p(u)\ptl_{p'}(v)-\ptl_{p'}(u)\ptl_p(v))+x_1^{\sgm_n}[(2-\ptl)(u)\ptl_n(v)-\ptl_n(u)(2-\ptl)(v)]\eqno(5.12)$$
for $u,v\in{\cal A}$ (cf. (4.22)). It can be verified that $({\cal A},[\cdot,\cdot]_K)$ forms a Lie algebra.
We call $({\cal A},[\cdot,\cdot]_K)$ a {\it  generalized Lie algebra of Contact type}.

\psp

{\bf Theorem 5.1}. {\it The pair} $({\cal A},[\cdot,\cdot]_K)$ {\it forms a simple Lie algebra}.
\psp

{\it Proof}. Let $I$ be a nonzero ideal of ${\cal A}$. Set
$${\cal A}'=\mbox{span}\:\{x^{\al,\vec{i}}\mid (\al,\vec{i})\in \G_1\times \vec{\cal J},\;i_n=0\}.\eqno(5.13)$$
In fact,
$${\cal A}'=\{u\in{\cal A}\mid\ptl_n(u)=0\}.\eqno(5.14)$$
\psp

{\it Step 1}. $I\bigcap {\cal A}'\neq\{0\}$.
\psp

We set
$${\cal A}_{[k]}=\mbox{span}\:\{x^{\al,\vec{i}}\mid (\al,\vec{i})\in\G\times\vec{\cal J},\;i_n\leq k\}\qquad\for\;\;k\in\Bbb{N};\qquad{\cal A}_{[-1]}=\emptyset.\eqno(5.15)$$
Define
$$\hat{k}=\min \{k\mid {\cal A}_{[k]}\bigcap I\neq\{0\}\}.\eqno(5.16)$$
Let $0\neq u\in {\cal A}_{[\hat{k}]}\bigcap I$. We write
$$u=u_0+u'\;\;\;\mbox{with}\;\;u'\in{\cal A}_{[\hat{k}-1]}\eqno(5.17)$$
and 
$$u_0=\sum_{(\al,\vec{i})\in\G\times\vec{\cal J},i_n=\hat{k}}a_{\al,\vec{i}}x^{\al,\vec{i}},\qquad a_{\al,\vec{i}}\in \Bbb{F}.\eqno(5.18)$$
Moreover, we define
$$\hat{\imath}(u)=|\{\al\mid a_{\al,\vec{i}}\neq 0\}|.\eqno(5.19)$$
Let
$$\hat{\imath}=\min\{\hat{\imath}(v)\mid 0\neq v\in {\cal A}_{[\hat{k}]}\bigcap I\}.\eqno(5.20)$$
Let $0\neq u\in {\cal A}_{[\hat{k}]}\bigcap I$ such that $\hat{\imath}(u)=\hat{\imath}$. Write $u$ as (5.17) and (5.18). 
\psp

{\it Case 1}. There exist $a_{\al,\vec{i}}\neq 0$ and $a_{\be,\vec{\cal J}}\neq 0$ such that $\al\in\G_1$ and $\be\not\in\G_1$.

Note that
$$[1,x^{\gm,\vec{l}}]_K=2\vf_n(\gm)x^{\gm+\sgm_n,\vec{l}}+2l_nx^{\gm+\sgm_n,\vec{l}-1_{[n]}}\eqno(5.21)$$
for $(\gm,\vec{l})\in\G\times\vec{\cal J}$ by  (5.11) and (5.12). In particular, we have
$$[1,x^{\al,\vec{i}}]_K\equiv 0,\;\;[1,x^{\be,\vec{\cal J}}]_K\equiv 2\vf_n(\be)x^{\be+\sgm_n,\vec{\cal J}}\;\;(\mbox{mod}\;{\cal A}_{[\hat{k}-1]}).\eqno(5.22)$$
Hence we have
$$0\neq [1,u]_K \in {\cal A}_{[\hat{k}]}\bigcap I,\;\;0<\hat{\imath}([1,u]_K)<\hat{\imath}(u)=\hat{\imath},\eqno(5.23)$$
which contradicts (5.20).
\psp

{\it Case 2}. $\al\in\G_1$ whenever $a_{\al,\vec{i}}\neq 0$.

If $\hat{k}=0$, then $u\in {\cal A}'$. So $I\bigcap {\cal A}'\neq \{0\}$. Assume that $\hat{k}>0$. By (5.21),
$0\neq [1,u]_K\in {\cal A}_{[\hat{k}-1]}\bigcap I$, which contradicts (5.16).
\psp

{\it Case 3}. There exists $a_{\al,\vec{i}}\neq 0$  such that  
$$(\vt(\al,\vec{i})+2)\vf_n(\al)+\vt(\al,\vec{i})\vf_n(\sgm_n)\neq 0,\eqno(5.24)$$ 
where
$$\vt(\al,\vec{i})=2-\sum_{p\in\mho_1}\vf_p(\al)-\sum_{q\in\mho_2}i_q.\eqno(5.25)$$

 Note that
$$[x_1^{\kappa},x^{\gm,\vec{l}}]_K=(2\vf_n(\gm)-\vf_n(\kappa)\vt(\gm,\vec{l}))x^{\gm+\kappa+\sgm_n,\vec{l}}+2l_nx^{\gm+\kappa+\sgm_n,\vec{l}-1_{[n]}}\eqno(5.26)$$
for $\kappa\in\G_2,\;\gm\in\G$ and $\vec{l}\in\vec{\cal J}$. Moreover, we write $\al=\al_1+\al_2$ with $\al_r\in \G_r$ by (5.2) and have 
$$0\neq [x_1^{-\al_2-\sgm_n},u]\in {\cal A}_{[\hat{k}]}\bigcap I,\;\;0<\hat{\imath}([x_1^{-\al_2-\sgm_n},u]_K)\leq \hat{\imath}\eqno(5.27)$$
and $[x_1^{-\al_2-\sgm_n},u]$ contains the following term:
$$a_{\al,\vec{i}}((\vt(\al,\vec{i})+2)\vf_n(\al)+\vt(\al,\vec{i})\vf_n(\sgm_n))x^{\al_1,\vec{i}}.\eqno(5.28)$$
Replacing $u$ by $[x_1^{-\al_2-\sgm_n},u]$, we go back to Cases 1 and 2.
\psp

{\it Case 4}. 
$$\al\not\in\G_1,\;\;(\vt(\al,\vec{i})+2)\vf_n(\al)+\vt(\al,\vec{i})\vf_n(\sgm_n)= 0\;\;\mbox{whenever}\;\;a_{\al,\vec{i}}\neq 0.\eqno(5.29)$$

In this case, $\G_2\neq\{0\}$. 
\psp

{\it Subcase 1}. There exists $a_{\al,\vec{i}}\neq 0$ such that $\vt(\al,\vec{i})\neq -2$.

We can choose $-\sgm_n\neq\gm\in\G_2$ such that
$$2\vf_n(\al)\neq \vf_n(\gm) \vt(\al,\vec{i})\eqno(5.30)$$
because $|\G_2|=\infty$. 
Then $[x_1^{\gm},u]$ contains the following term
$$a_{\al,\vec{i}}(2\vf_n(\al)-\vf_n(\gm) \vt(\al,\vec{i}))x^{\al+\gm+\sgm_n,\vec{i}}.\eqno(5.31)$$
Thus, replacing $u$ by $[x_1^{\gm},u]$, we go back to Case 3 because $\vf_n|_{\G_2}$ is injective.
\psp

{\it Subcase 2}. $\vt(\al,\vec{i})= -2$ whenever $a_{\al,\vec{i}}\neq 0$.

In this subcase, $\sgm_n=0$ by (5.29). Moreover, (5.21) becomes
$$[1,x^{\gm,\vec{l}}]_K=2\ptl_{n}(x^{\gm,\vec{l}})=2\vf_n(\gm)x^{\gm,\vec{l}}+2l_{n}x^{\gm,\vec{l}-1_{[n]}}\eqno(5.32)$$
for $(\gm,\vec{l})\in\G\times\vec{\cal J}$.
For any $\be\in\G_2$, we let
$${\cal A}_{[\be]}=\mbox{span}\:\{x^{\al+\be,\vec{i}}\mid (\al,\vec{i})\in \G_1\times \vec{\cal J}\}.\eqno(5.33)$$
Then we have:
$${\cal A}=\bigoplus_{\be\in\G_2}{\cal A}_{[\be]},\;\;\;{\cal A}_{[\be]}=\{u\in {\cal A}\mid (\ad_1-2\vf_n(\be))^j(u)=0\;\mbox{for some}\;j\in\Bbb{N}\}\eqno(5.34)$$
by (5.32) and (5.33). 

Let $a_{\al,\vec{i}}\neq 0$ be such that $|\vec{i}|$ is maximal (cf. (4.14)). Again we write $\al=\al_1+\al_2$ with $\al_r\in\G_r$. 
If ${\cal J}_1=\Bbb{N}$, then we have
$$[x^{\al,\vec{i}},x^{-\al_2,1_{[1]}}]_K= \vf_n(\al)x^{\al_1,\vec{i}+1_{[1]}}-\vf_{m+1}(\al)x^{\al_1,\vec{i}}-i_1x^{\al,\vec{i}-1_{[1]}}\eqno(5.35)$$
by (5.11), (5.12) and (5.29). Thus $[u,x^{-\al_2,1_{[1]}}]$ contains the following term:
$$a_{\al,\vec{i}}\vf_n(\al)x^{\al_1,\vec{i}+1_{[1]}}\in {\cal A}_{0_{\G_2}}.\eqno(5.36)$$
Since $[x^{-\al_2,1_{[1]}},u]_K\in I$, (5.34) and (5.36) imply
$$I\bigcap {\cal A}_{0_{\G_2}}\neq \{0\}.\eqno(5.37)$$
Pick any $0\neq v\in I\bigcap {\cal A}_{0_{\G_2}}$. Applying $\ad_1$ on $v$ repeatedly, we can get a nonzero element $w\in {\cal A}'$ by (5.32). 

If ${\cal J}_1=\{0\}$, then $\vf_1\not\equiv 0$. So $\sgm_1\neq 0$ by the assumption in (5.6). Note
$$[x_1^{\sgm_1-\al_2},x^{\al,\vec{i}}]_K\equiv(\vf_1(\al)-\vf_{m+1}(\al))x^{\al_1+2\sgm_1,\vec{i}}+\vf_n(\al)x^{\al_1+\sgm_1,\vec{i}}\;\;(\mbox{mod}\;\;{\cal A}_{[\hat{k}-1]})\eqno(5.38)$$
by (5.8), (5.11), (5.12) and (5.29). Thus $[x_1^{\sgm_1-\al_2},u]_K$ contains the following term
$$a_{\al,\vec{i}}\vf_n(\al)x^{\al_1+\sgm_1,\vec{i}}\in {\cal A}_{0_{\G_2}}.\eqno(5.39)$$
Since $[x_1^{\sgm_1-\al_2},u]_K\in I$, (5.34) and (5.39) imply (5.37).
Thus $I\bigcap {\cal A}'\neq\{0\}$.
\psp

{\it Step 2}. $1\in I$.
\psp

Note that for $u,v\in {\cal A}'$, we have
$$[u,v]_K=\sum_{p=1}^mx_1^{\sgm_p}[\ptl_p(u)\ptl_{m+p}(v)-\ptl_{m+p}(u)\ptl_p(v)]\eqno(5.40)$$
by (5.11) and (5.12). Moreover, (5.40) is a special case of (4.12) with $m_1=m,\;\phi\equiv 0$ and ${\cal A}$ replaced by ${\cal A}'$. By the proof of Theorem 4.2, we have
$$I\bigcap {\cal A}'=\Bbb{F}\eqno(5.41)$$
or
$$ I\bigcap {\cal A}'={\cal A}'\qquad\mbox{if}\;\;{\cal J}_p=\Bbb{N}\;\;\mbox{for some}\;\;p\in\ol{1,2m}\eqno(5.42)$$
or
$$I\bigcap {\cal A}'=\mbox{span}\;\{x_1^{\al}\mid \sum_{p=1}^m\sgm_p\neq\al\in\G_1\}\;\;\mbox{if}\;\;{\cal J}_p=\{0\}\;\;\mbox{for any}\;\;p\in\ol{1,2m}\eqno(5.43)$$
(cf. (4.70), (4.79)). Note by (5.1), (5.3) and (5.6),
$$\sum_{p=1}^m\sgm_p\neq 0\eqno(5.44)$$
in the third case. Thus we always have $1\in I$ in all the three cases (5.41)-(5.43).
\psp

{\it Step 3}. $I={\cal A}$ from $1\in I$.
\psp

By (5.21) and induction on $l_n$, we can prove $I={\cal A}$ if ${\cal J}_n=\Bbb{N}$.
Assume ${\cal J}_n=\{0\}$. Then $\vf_n\not\equiv 0$ by the assumptions in (5.1). 
Note that (5.21) shows
$$x^{\al_1+\al_2,\vec{i}}\in I\qquad\for\;\al_1\in\G_1,\;\sgm_n\neq \al_2\in\G_2,\;\vec{i}\in\vec{\cal J}.\eqno(5.45)$$
For any $0,\sgm_n\neq \kappa\in \G_2$ and $(\be,\vec{\cal J})\in \G_1\times\vec{\cal J}$ such that $\vt(\be,\vec{\cal J})\neq -2$ , we have
$$x^{\be+\sgm_n,\vec{\cal J}}=[\vf_n(\kappa)(2+\vt(\be,\vec{\cal J}))]^{-1}[x^{-\kappa}_1,x^{\be+\kappa,\vec{\cal J}}]_K\in I.\eqno(5.46)$$
Since  $\vt(0,1_{[1]})=1$ if ${\cal J}_1=\Bbb{N}$ and $\vt(\sgm_1,\vec{0})=3$ if $\vf_1\not\equiv 0$ by (5.8) and (5.25), (5.45) and (5.46) imply  either $x_2^{1_{[1]}}\in I$ or $x_1^{\sgm_1}\in I$. Moreover, the arguments in Cases 1 and 2 of Step 2 in the proof of Theorem 4.2 show that $I={\cal A}$ if $\vec{\cal J}\neq\{\vec{0}\}$ and  
$$x_1^{\al_1+\al_2}\in I\qquad\for\;\;\sum_{p=1}^m\sgm_p\neq \al_1\in\G_1,\;\al_2\in\G_2\eqno(5.47)$$
if $\vec{\cal J}=\{\vec{0}\}$. When $\vec{\cal J}=\{\vec{0}\}$, (5.45), (5.46) and (5.47) imply $I={\cal A}$ because $\vt(\sgm,\vec{0})=2(m+1)\neq -2$ by (5.8) and (5.25). $\qquad\Box$
\psp

{\bf Example}.  Suppose that we have picked (2.5). Let $k$ be the number of ${\cal J}_p =\{0\}$ with $p\in\ol{1,n}$. Pick an integer $\ell$ such that $k\leq \ell\leq n$. We define
$$\zeta_p(\al_1,...,\al_{\ell})=\al_p\qquad\for\;\;p\in\ol{1,\ell},\;(\al_1,\al_2,...,\al_{\ell})\in\Bbb{F}^{\ell}\eqno(5.48)$$
and
$$\zeta_q\equiv 0 \qquad\for\;\;q\in\ol{\ell+1,n}.\eqno(5.49)$$
Pick any permutation $\iota$ on $\ol{1,n}$ such that
$$\iota(p)\leq\ell\;\;\mbox{if}\;\;{\cal J}_p=\{0\}\;\;\mbox{with}\;\;p\in\ol{1,n}.\eqno(5.50)$$
Moreover, we can assume $\iota(n)=\ell$ if $\iota(n)\leq \ell$. 
Take $\G$ to be an additive subgroup of $\Bbb{F}^{\ell}$ such that $\G\supset\Bbb{Z}^{\ell}$ and
$$\G\supset \{(0,...,0,\al_{\ell})\mid (\al_1,...,\al_{\ell})\in \G\}\;\;\mbox{if}\;\;\iota(n)=\ell.\eqno(5.51)$$
We let
$$\vf_p=\zeta_{\iota(p)}|_{\G}\qquad\for\;\;p\in\ol{1,n}.\eqno(5.52)$$
Then (5.1)-(5.3) and (5.5) hold. 

In particular, we can take (2.34) and (2.35).

\vspace{1cm}

\noindent{\Large \bf References}

\hspace{0.5cm}

\begin{description}

\item[{[DZ]}] D. \v{Z}. Dokovick and K. Zhao, Generalized Cartan type W Lie algebras in characteristc 0, {\it J. Algebra} {\bf 195} (1997), 170-210.

\item[{[B]}] R. Block, On torsion-free abelian groups and Lie algebras, {\it Proc. Amer. Math. Soc.} {\bf 9} (1958), 613-620.

\item[{[Ka1]}] V. G. Kac, A description of filtered Lie algebras whose associated graded Lie algebras are of Cartan types, {\it Math. of USSR-Izvestijia} {\bf 8} (1974), 801-835.

\item[{[Ka2]}] V. G. Kac, Lie superalgebras, {\it Adv. Math.} {\bf 26} (1977), 8-96.

\item[{[Ka3]}] V. G. Kac, {\it Vertex algebras for beginners}, University lectures series, Vol {\bf 10}, AMS. Providence RI, 1996.

\item[{[Ka4]}] V. G. Kac, Superconformal algebras and transitive group actions on quadrics, {\it Commun. Math. Phys.} {\bf 186} (1997), 233-252.

\item[{[KL]}]  V. G. Kac and J. W. Leur, On classification of superconformal algebras, in S. J. et al. eds. {\it String 88}, World Sci. (1989), 77-106.

\item[{[KT]}] V. G. Kac and I. T. Todorov, Superconformal current algebras and their unitary representations, {\it Commun. Math. Phys}. {\bf 102} (1985), 337-347.

\item[{[K]}] N. Kawamoto, Generalizations of Witt algebras over a field of characteristic zero,
{\it Hiroshima Math. J.} {\bf 16} (1986), 417-462.

\item[{[O1]}]
J. Marshall Osborn, Infinite dimensional Novikov algebras of characteristic 0, {\it J. Algebra} {\bf 167} (1994), 146-167.

\item[{[O2]}] J. Marshall Osborn, New simple infinite-dimensional Lie algebras of characteristic 0, {\it J. Algebra} {\bf 185} (1996), 820-835.

\item[{[P]}] D. P. Passman, Simple Lie algebras of Witt type, {\it J. Algebra} {\bf 206} (1998), 682-692.

\item[{[SF]}] H. Strade and R. Farnsteiner, {\it Modular Lie Algebras and Their Representations,} Marcel Derkker, Inc., 1988.

\item[{[X1]}] X. Xu, On simple Novikov algebras and their irreducible modules, {\it J. Algebra} {\bf 185} (1996), 905-934.

\item[{[X2]}] X. Xu, Variational calculus of supervariables and related algebraic structures, {\it submitted,} preprint was circulated in January 1995.

\item[{[X3]}] X. Xu, Novikov-Poisson algebras, {\it J. Algebra} {\bf 190} (1997), 253-279.

\item[{[X4]}] X. Xu, Quadratic conformal superalgebras, {\it submitted}.

\item[{[X5]}] X. Xu, {\it Introduction to Vertex Operator Superalgebras and Their Modules}, Kluwer Academic Publishers, 1998.

\item[{[X6]}] X. Xu, Generalizations of Block algebras, {\it submitted}.

\end{description}
\end{document}